\documentclass[12pt]{article}
\usepackage{amsfonts}
\usepackage{epsfig}
\title{Unbounded Orbits in the Plaid Model}
\author{Richard Evan Schwartz \thanks{\hskip 5 pt Supported by 
N.S.F. Research Grant DMS-1204471}}

\newtheorem{theorem}{Theorem}[section]

\newtheorem{lemma}[theorem]{Lemma}

\newtheorem{corollary}[theorem]{Corollary}

\def\startproof{{\bf {\medskip}{\noindent}Proof: }}

\def\endproof{$\spadesuit$  \newline}

\def\H{\mbox{\boldmath{$H$}}}%
\def\N{\mbox{\boldmath{$N$}}}%
\def\Q{\mbox{\boldmath{$Q$}}}%
\def\R{\mbox{\boldmath{$R$}}}%
\def\Z{\mbox{\boldmath{$Z$}}}%

\begin{document}
\maketitle

\begin{abstract}
This is a sequel to [{\bf S0\/}].
In this paper I prove that the plaid
model has unbounded orbits at all
irrational parameters.  This result
is closely related to my result [{\bf S1\/}] that
outer billiards has unbounded orbits
with respect to any irrational kite. 
\end{abstract}

\section{Introduction}

I introduced the plaid model in
[{\bf S0\/}], and this paper is a sequel.
The plaid model is closely
related to outer billiards on kites
[{\bf S1\/}], somewhat
related to DeBruijn's pentagrids [{\bf DeB\/}],
and also somewhat related to
corner percolation and P. Hooper's
Truchet tile system [{\bf H\/}].

The plaid model has $2$ descriptions.
For a rational parameter
$p/q \in (0,1)$ with $pq$ even,
there is one description in terms
of intersection points on grids
of lines.  See \S 2.  The output
of the plaid model is a union of
embedded polygons in the plane
which we call {\it plaid polygons\/}.

At the same time, there is a
$3$ dimensional polyhedron exchange
transformation 
$\widehat X_P$ which exists for
all $P \in (0,1)$.
  This PET has
the property that
the plaid polygons associated to
the parameter $p/q$ correspond
to the so-called 
vector dynamics of a
distinguished set of orbits of
$\widehat X_P$ for $P=2p/(p+q)$.
This is [{\bf S0\/}, Theorem 1.1].
See \S 3.

We can interpret the vector
dynamics at irrational parameters
as (Hausdorff) limits of these
plaid polygons.  We call these
limits {\it plaid paths\/}.
In other words, the plaid paths
are the vector dynamics for the PET
at all parameters, and for even
rational parameters they are
embedded closed loops.
All this is explained in \S 2-3.

The purpose of this paper is to prove
\begin{theorem}[Unbounded Orbits]
For every irrational $P \in (0,1)$ the
PET $\widehat X_P$ has an infinite orbit.
The corresponding plaid path has infinite
diameter projections in both coordinate
directions.
\end{theorem}
The reason one might want the coordinate
projections of a plaid path to be
unbounded is that the projection into
the $x$-coordinate has a dynamical
interpretation in terms of outer billiards
orbits.  According to the Quasi-Isomorphism
Conjecture in [{\bf S0\/}], a plaid path with
unbounded $x$-projection corresponds to
an unbounded orbit for outer billiards
on the kite with the same parameter.

In addition to the Unbounded Orbits Theorem,
we will prove a number of results about
how the set of plaid polygons associated
to a given parameter contains pieces of
the plaid polygons associated to a different
and closely related parameter.
I mean {\it closely related\/} in the
sense of Diophantine approximation.
These copying results are
difficult to state without a buildup
of terminology; in the paper they are
called the Copy Theorem, and the
three Copy Lemmas.

The proof of the Unbounded Orbits Theorem here
is similar in spirit to the proof in [{\bf S1\/}] of the
Unbounded Orbits Theorem for outer billiards
on kites.
One difference is that the proof here is
traditional whereas the proof in [{\bf S1\/}]
was heavily computer-assisted. (The proof here
quotes some results from [{\bf S0\/}], but
these results are also proved in a traditional way.)
The plaid model is somehow more straightforward
and transparent than the corresponding object
from [{\bf S1\/}], namely the
{\it arithmetic graph\/}, and so it is
possible to do all the calculations by hand.

On the other hand, the results here are 
like the results in [{\bf S1\/}] in that we
discovered them through experimentation, and
checked practically every step on many
examples before trying for proofs. As
I remarked in [{\bf S1\/}], the massive
reliance on the computer for motivation
is a virtue rather than a bug: I checked
all the identities in the paper on
thousands of examples.

The Quasi-Isomorphism Conjecture is a piece of
unfinished business from
[{\bf S0\/}]. This result says that,
for every parameter, there is a bijection
between the set of plaid polygons and the
set of components of the arithmetic graph
which differs from the identity map by
at most $2$ units.
I hope to prove the
Quasi-Isomorphism Theorem next.  Given the
truth of the Quasi-Isomorphism Conjecture, the
Unbounded Orbits Theorem here is very nearly
equivalent to the one in [{\bf S1\/}], but
perhaps there are a few small details which
would need to be worked out. 

Here is an overview of the paper.
\begin{itemize}
\item In \S 2, we discuss the definition of the
plaid model in terms of grids of lines.  In
particular, we give a new interpretation
of the plaid model in terms of triple
intersections of lines which are decorated
with signs.
\item In \S 3, we discuss the plaid model in terms
of polyhedron exchange transformations.
The material here is the same as in
[{\bf S0\/}] except that we discuss
geometric limits in more detail.
\item In \S 4 we prove the Unbounded
Orbits Theorem modulo some elementary
number theory, and and two auxiliary
results, the Box Lemma and the Copy Theorem.
\item In \S 5 we take care of the
number-theoretic component of the proof.
In particular, we prove Lemma \ref{omnibus},
a multi-part result which summarizes a
number of small number-theoretic
identities and inequalities between
rational numbers.  
\item In \S 6 we reduce the Box Lemma and
the Copy Theorem
to three similar but smaller results,
which we call the Weak, Strong, and Core
Copy Lemmas.
\item In \S 7 we establish the Weak and Strong
Copy Lemmas.  This is an application of the
signed crossing interpretation of the
Plaid Model. 
\item In \S 8 we establish the Core Copy Lemma.
This is another application of the signed
crossing interpretation. 
\end{itemize}

\newpage

\section{The Grid Description}

\subsection{Basic Definitions}

\noindent
{\bf Even Rational Parameters:\/}
We will work with 
$p/q \in (0,1)$ with $pq$ even.
We call such numbers {\it even rational parameters\/}.
We will work closely with the auxiliary quantities
\begin{equation}
\omega=p+q, \hskip 30 pt
P=\frac{2p}{\omega}, \hskip 30 pt
Q=\frac{2q}{\omega}.
\end{equation}
\newline
\newline
{\bf Four Families of Lines\/}
We consider $4$ infinite families of lines.
\begin{itemize}
\item $\cal H$ consists of horizontal lines having integer $y$-coordinate.
\item $\cal V$ consists of vertical lines having integer $x$-coordinate.
\item $\cal P$ is the set of lines of slope $-P$ having
integer $y$-intersept.
\item $\cal Q$ is the set of lines of slope $-Q$ having
integer $y$-intersept.
\end{itemize}

\noindent
{\bf Adapted Functions:\/}
Let $\Z_0$ and $\Z_1$ denote the
sets of even and odd integers respectively.
We define the following $4$ functions:
\begin{itemize}
\item $F_H(x,y)=2Py$ mod $2\Z$
\item $F_V(x,y)=2Px$ mod $2\Z$.
\item $F_P(x,y)=Py+P^2x+1$ mod 2$\Z$
\item $F_Q(x,y)=Py+PQx+1$ mod 2$\Z$.
\end{itemize}
We call $F_H$ and $F_V$ {\it capacity functions\/}
and $F_P$ and $F_Q$ {\it mass functions\/}.
We shall never be interested in the
inverse images of $\Z_1$, so
we will always normalize so that our
functions take values in $(-1,1)$.
\newline
\newline
{\bf Intersection Points:\/}
Recall that $\omega=p+q$.  Let
\begin{equation}
\Omega_j=\Z_j \omega^{-1}.
\end{equation}
We observe that, for each
index $A \in \{H,V\}$ and
$B \in \{P,Q\}$, we have
$F_A^{-1}(\Omega_0)=\cal A$ and
$F_B^{-1}(\Omega_1)=\cal B$.

We call $z$ an {\it intersection point\/}
if $z$ lies on both an $A$ line and a $B$ line.
In this case,
we call $z$ {\it light\/} if
\begin{equation}
\label{LIGHT}
|F_B(z)|<|F_A(z)|, \hskip 30 pt
F_A(z)F_B(z)>0.
\end{equation}
Otherwise we call $z$ {\it dark\/}.
We say that $B$ is the {\it type\/}
of $z$.  Thus $z$ can be light or dark,
and can have type P or type Q.  When
$z$ lies on both a $\cal P$ line and
a $\cal Q$ line, we say that $z$
has both types.
\newline
\newline
{\bf Capacity and Mass:\/}
For any index $C \in \{V,H,P,Q\}$, we
assign two invariants to a line
$L$ in $\cal C$, namely
$|(p+q)F_C(L)|$ and the sign of
$F_C(L)$.  We call the second invariant
the {\it sign\/} of $L$ in all cases.
 When $A\in\{V,H\}$, we call
the first invariant the {\it capacity\/} of $L$.
When $A \in \{P,Q\}$ we call the
first invariant the {\it mass\/} of $L$.
The masses are all odd integers in
$[1,p+q]$ and the capacities are
all even integers in $[0,p+q-1]$.
Thus, an intersection point is light
if and only if the intersecting
lines have the same sign, and the
mass of the one line is less than
the capacity of the other.
\newline
\newline
{\bf Double Counting Midpoints:\/}
We introduce the technical rule that
we count a light point twice if it
appears as the midpoint of a horizontal
unit segment.  The justification for this
convention is that such a point is 
always a triple intersection of a
$\cal P$ line, a $\cal Q$ line, and
a $\cal H$ line; and moreover this point
would be considered light either when
computed either with respect to the
$\cal P$ line or with respect to the
$\cal Q$ line.  See Lemma \ref{weak} below.
\newline
\newline
{\bf Good Segments:\/}
We call the edges of the unit squares
{\it unit segments\/}.  Such segments,
of course, either lie on $\cal H$
lines or on $\cal V$ lines.  We say
that a unit segment is {\it good\/} if it
contains exactly one light point.
  We say that a unit square
is {\it coherent\/} if it contains either
$0$ or $2$ good segments.  We say that the
plaid model is {\it coherent\/} at if all
squares are coherent for all parameters.
Here is the fundamental theorem concerning
the plaid model.  [{\bf S0\/}, Fundamental Theorem] 
says that the plaid model is coherent at all parameters.
\newline
\newline
\noindent
{\bf Plaid Polygons:\/}
In each unit square we draw the line segment which connects
the center of the square with the centers of
its good edges.   The squares with no good edges
simply remain empty. We call these polygons
the {\it plaid polygons\/}.

\subsection{Symmetries}
\label{plaidsymm}

\noindent
{\bf The Symmetry Lattice\/}
We fix some even rational parameter $p/q$.
Let $\omega=p+q$ as above.
Let $L \subset \Z^2$ denote the
lattice generated by the two vectors
\begin{equation}
(\omega^2,0), \hskip 30 pt
(0,\omega).
\end{equation}
We call $L$ the {\it symmetry lattice\/}.
\newline
\newline
{\bf Blocks:\/}
We define the square $[0,\omega]^2$
to be the {\it first block\/}.
The pictures above always show the first block.
In general, we define a {\it block\/} to be
a set of the form $B_0 + \ell$, where $B_0$ is the
first block and $\ell \in L$.
With this definition, the lattice $L$
permutes the blocks.  
We define the {\it fundamental blocks\/} to
be $B_0,...,B_{\omega-1}$, where
$B_0$ is the first block and
\begin{equation}
B_k=B_0+(k\omega,0).
\end{equation}
The union of the fundamental blocks is a
fundamental domain for the action of $L$.
We call this union the {\it fundamental domain\/}.
\newline
\newline
{\bf Translation Symmetry:\/}
All our functions $F$ are $L$-invariant.
That is, $F(v+\ell)=F(v)$ for all $\ell \in L$.
This follows immediately from the definitions.
Note that a unit segment in the
boundary of a block does not have any
light points associated to it
because $F_H=0$ on the horizontal
edges of the boundary and $F_V=0$ on
the vertical edges of the boundary. 
(A vertex of such a segment might be a
light point associated to another
unit segment incident to it but this
doesn't bother us.)
Hence, every plaid polygon is
contained in a block, and is translation equivalent
to one in a fundamental block.
\newline
\newline
{\bf Reflection Symmetry:\/}
Let $\widehat L$ denote the group
of isometries generated by
reflections in the horizontal and
vertical midlines (i.e. bisectors)
of the blocks.  In \S 2 we proved
that the set of plaid polygons
is invariant under the action
of $\widehat L$.
Regarding the types of the intersection
point (P or Q), the reflections in $\widehat L$
have the following action:
\begin{itemize}
\item They preserve the types of the intersection
points associated to the horizontal lines.

\item They interchange the types of the intersection
points associated to the vertical lines.
\end{itemize}

\subsection{An Alternate Description}

We have defined the plaid model in terms of
the $\cal H$ and $\cal V$ lines, and the
$\cal P$ and $\cal Q$ lines.  Our description
favored lines of negative slope.
To clarify the situation, we refer to the
$\cal P$ and $\cal Q$ lines as the
${\cal P\/}_-$ and ${\cal Q\/}_-$ lines.
We define the ${\cal P\/}_+$ and
${\cal Q\/}_+$ lines to be the lines 
having integer $y$-intercept and
slope $+P$ and $+Q$ respectively

We define the mass of a ${\cal P\/}_+$ line $L_+$
so that it equals the mass of the
${\cal P\/}_-$ line $L_-$ having the same $y$-intercept,
and we define the sign of the $L_+$ to be
opposite the sign of $L_-$.
We make the same definitions for the
${\cal Q\/}_+$ lines.   

Here is a
result that is implicitly contained in
[{\bf S0\/}], but we want to make it explicit here.

\begin{lemma}[Vertical]
Let $z$ be an intersection point on a $\cal V$ line.
Then either $z$ is the intersection of a
${\cal P\/}_+$ line with a ${\cal Q\/}_-$ line or
$z$ is the intersection of a
${\cal P\/}_-$ line with a ${\cal Q\/}_+$ line.
Moreover $z$ is a light particle if and only if
the three lines containing $z$ have the same sign.
\end{lemma}

\startproof
We have the functions
$$F_{P,-}(x,y)=[Py+P^2x+1]_2, \hskip 30 pt
  F_{Q,-}(x,y)=[Py+PQx+1]_2,$$
corresponding to the ${\cal P\/}_-$ and
${\cal Q\/}_-$ lines.  There are similar
functions corresponding to the
${\cal P\/}_+$ and ${\cal Q\/}_+$ lines,
namely:
$$F_{P,+}(x,y)=[-Py+P^2x+1]_2, \hskip 30 pt
  F_{Q,+}(x,y)=[-Py+PQx+1]_2,$$
Let's say that $z$ lies on the intersection
of a $\cal V$ line and a ${\cal P\/}_-$ line.
Then $\omega F_{P,-}(z)$ is an odd integer
in $(-\omega,\omega)$ and
$\omega F_{V,-}(z)$ is an even integer in
$(-\omega,\omega)$.

But then we have the identity
\begin{equation}
\label{key}
F_{P,-}+F_{Q,+}=[2Px]_2=F_H.
\end{equation}
This tells us that $\omega F_{Q,+}(z)$ is also an
odd integer and so $z$ lies on some ${\cal Q\/}_+$
line.  If all three lines have the same sign, then
Equation \ref{key} forces $F_{P,-}(z)<F_H(z)$.
So $z$ is a light particle.
If the ${\cal P\/}_-$ line and the $\cal V$ line
have opposite signs then by definition $z$ is dark.
If the ${\cal Q\/}_+$ line and the $\cal V$ line
have opposite signs, then Equation \ref{key}
forces $F_{P,-}(z)>F_H(z)$ and again $z$ is dark.
\endproof

\noindent
{\bf The Inert Principle:\/}
Equation \ref{key} has one more consequence of
importance to us.  The $\cal P$ and $\cal Q$ lines
through the points in $\omega\Z$ have mass
$\omega$, but it is impossible to give them 
a sign. We call these lines {\it inert\/}
and we call the intersection points on
them inert as well.
All the inert intersection points are dark.
Thanks to Equation \ref{key}, one of
the following is true for each inert
intersection point $z$.
\begin{itemize}
\item $z$ lies in the boundary of
a block, and the two slanting lines
containing it are both inert.
\item The other slanting line
containing $z$ has the opposite
sign as the horizontal/vertical
line containing $z$.
\end{itemize}
We call this dichotomy the
{\it Inert Principle\/}.

The Vertical Lemma and the Inert Principle tell
us that the types (light versus dark) of the
vertical intersections in the plaid model are
entirely determined by the signs of the grid
lines.  The same goes for the horizontal
intersections, though the condition on the
signs is different:  The ${\cal H\/}/{\cal V\/}$
lines and the $(-)$ lines should have the same
sign and the $(+)$ lines should have the
opposite sign.  Of course, we have set our
sign conventions to favor a simpler description
of the vertical intersections.

\subsection{The Mass and Capacity Sequences}

Now we will package the information in the
last section in a way that will be more useful to us.
Let

\begin{equation}
\Sigma=[x_0,x_1] \times [y_0,y_1], \hskip 30 pt
\end{equation}
We assume that $\Sigma$ does not intersect
the vertical boundary of a block.  More
precisely, for this paper we will have
$1=x_0<x_1<\omega$.

All the $\cal V$ lines which intersect $\Sigma$ have
positive capacity. The signs could be
positive or negative.
We define two sequences 
$\{c_j\}$ and $\{m_j\}$ w.r.t. $\Sigma$.
We have the {\it capacity sequence\/}
\begin{equation}
c_j=[2Pj]_2, \hskip 20 pt
j=x_0,...,x_1.
\end{equation}
The terms of the capacity sequence all lie in
$(-1,1)$.

We also have the {\it mass sequence\/}
\begin{equation}
m_j=[P_j+1]_2, \hskip 30 pt
j=y_0-2x+1,...,y_1+2x-1,
\hskip 30 pt x=x_0-x_1.
\end{equation}

\noindent
{\bf Remarks:\/} \newline
(i)
Notice that the indices for the mass sequence
start below the bottom edge of $\Sigma$ (so to speak)
and end above it.  This is important for us
for reasons which will become clear momentarily.
\newline
(ii)
We will allow $j \in \omega\Z$ in the
mass sequence.  Such terms do not have
a well-defined sign. The corresponding
slanting lines are inert.  When we speak
of the signs of the terms of the mass sequence,
we mean to ignore these terms.
\newline
(iii) 
Really we only care about the signs in the mass
and capacity sequences, but for the purposes
of running certain kinds of arguments it is
useful to keep track of the numerical
values as well.
\newline

\begin{lemma}
The shade of any vertical intersection point
in $\Sigma$ is determined by the signs of the
terms in the mass and capacity sequences.
\end{lemma}

\startproof
Since the slanting lines in the plaid model
have slopes in $(-2,2)$, every slanting line
which contains a point of $\Sigma$ intersects
the $y$-axis in the interval
$$\{0\} \times [y_0-2x+1,y_1+2x-1].$$
Hence, every vertical intersection
point lies on $3$ lines whose signs are
all determined by the mass and capacity
sequence.  The Vertical Lemma and the
Inert Principle allow us to determine 
whether the intersection point is light
or dark.
\endproof

\begin{corollary}
\label{match00}
Suppose that we know how the plaid
polygons intersect a single $\cal H$
line inside $\Sigma$.  Then the
intersection of the plaid polygons
with $\Sigma$ is determined by the
signs of the mass and capacity sequences.
\end{corollary}

\startproof
We will suppose that we know how the plaid
polygons intersect the bottom edge of
$\Sigma$.  The case for any other edge
has a similar treatment.
Let $Q$ be some unit square in $\Sigma$
for which we have not yet determined the
plaid model inside $Q$.  We can take $Q$ to be
as low as possible.  But then we know
how the plaid model intersects the bottom
edge of $Q$, and the signs of the mass and
capacity sequences determine how the
tiling intersects the left and right edges.
But then the Fundamental Theorem for the
plaid model tells us that $\partial Q$
intersects the plaid polygons in either
$0$ or $2$ points.  This allows us to
determine how the plaid polygons
intersect the top edge of $Q$.
\endproof

\subsection{A Matching Criterion}
\label{MATCH}

We mean to define all the objects in
the previous section with respect to
parameter $p/q$ and $p'/q'$.
Our notation convention will be that
the object $X$ corresponds to $p/q$
whenever the same kind of object
$X'$ corresponds to $p'/q'$.

Let $\Pi$ denote the union of plaid
polygons with respect to $p/q$ and
let $\Pi'$ denote the union of
plaid polygons with respect to $p'/q'$.
Suppose that $\Sigma$ and $\Sigma'$ are
rectangles that are equivalent via a
vertical translation $\Upsilon$.
That is, $\Upsilon$ preserves the
$y$-axis and $\Upsilon(\Sigma')=\Sigma$.
To be concrete, say that
\begin{equation}
\Upsilon(x,y)=(x,y+\xi).
\end{equation}
We will give a criterion which guarantees
that 
\begin{equation}
\label{perfect}
\Upsilon(\Sigma' \cap \Pi')=\Sigma \cap \Pi.
\end{equation}

\noindent
{\bf Arithmetic Alignment:\/}
A necessary condition for Equation
\ref{perfect} is that the signs of
$\{c_j\}$ are the same as the
corresponding signs of $\{c'_j\}$
and the signs of $\{m_j\}$ are
the same as the corresponding
signs of $\{m'_j\}$.  More precisely,
\begin{equation}
\label{arithmetic}
{\rm sign\/}(c'_j)={\rm sign\/}(c_{j+\xi}), \hskip 30 pt
{\rm sign\/}(m'_j)={\rm sign\/}(m_{j+\xi}).
\end{equation}
We say that
$(\Sigma,\Pi)$ and $(\Sigma',\Pi')$ are
{\it arithmetically aligned\/} of
Equation \ref{arithmetic} holds
for all relevant indices.

It seems plausible that arithmetic
alignment is sufficient for Equation
\ref{perfect} to hold, but we don't
have a proof.  We need more ingredients
to make things work out cleanly.
\newline
\newline
{\bf Geometric Alignment:\/}
There is a natural correspondence between the
vertical intersection points in $\Sigma$ and
the vertical intersection points in $\Sigma'$.
Let $z'$ be a vertical intersection point in
$\Sigma$.  Let $\{i',j'\}$ be the pair of
indices so that the slanting lines 
through $(0,i')$ and $(0,j')$ contain $z'$.  We
let $z$ denote the intersection of the slanting
lines, of the same type,
through $(0,i)$ and $(0,j)$.  Here
$i=i'+\xi$ and $j=j'+\xi$.
We say that $z'$ and $z$ are
{\it geometrically aligned\/} of
$\Upsilon(z')$ and $z$ are contained in the
same unit vertical segment of $\Sigma$.
We say that $(\Sigma,\Pi)$ and
$(\Sigma',\Pi')$ are {\it geometrically aligned\/}
if $z$ and $z'$ are geometrically aligned for
every vertical intersection point $z' \in \Sigma'$.

It seems very likely that arithmetic and 
geometric alignment together imply
Equation \ref{perfect} but we don't have
a proof. We need one more small ingredient.
\newline
\newline
{\bf Weak Horizontal Alignment:\/}
We say that $(\Sigma,\Pi)$ and $(\Sigma',\Pi')$ are
{\it weakly horizontally aligned\/} if
there are $\cal H$ lines $H'$ and $H$ such that
\begin{equation}
\Upsilon(\Sigma \cap H \cap P)=\Sigma' \cap H' \cap \Pi'.
\end{equation}
In other words, the tilings look the same on a
single horizontal segment.  This is exactly the
criterion which appears in
Corollary \ref{match00}.

Now we come to our Matching Criterion.
\begin{lemma}[Matching Criterion]
Suppose that
\begin{itemize}
\item $(\Sigma,\Pi)$ and $(\Sigma',\Pi')$ are
weakly horizontally aligned.
\item $(\Sigma,\Pi)$ and $(\Sigma',\Pi')$ are
geometrically aligned.
\item $(\Sigma,\Pi)$ and $(\Sigma',\Pi')$ are
arithmetically aligned.
\end{itemize}
Then $\Upsilon(\Sigma' \cap \Pi')=\Sigma \cap \Pi$.
\end{lemma}

\startproof
Given Corollary \ref{match00}, this is
practically a tautology.
The procedure given in
Corollary \ref{match00} assigns exactly the
same tiles to $\Sigma \cap \Pi$ as it
does to $\Upsilon(\Sigma' \cap \Pi')$.
\endproof

\subsection{The Big Polygon}
\label{tune}

We define $\tau$ to be the unique
integer solution in $(0,\omega/2)$ to
\begin{equation}
2p\tau \equiv \pm 1 \hskip 10 pt
{\rm mod\/} \hskip 5 pt \omega.
\end{equation}
We call $\tau$ the {\it tune\/} of
the parameter $p/q$.

\begin{lemma}
\label{anchor}
For $k=0,...,(\omega-1)/2$, the
lines of capacity $2k$ have the form
$x=\pm k \tau$ and
$y=\pm k \tau$.
For $k=1,3,...,(\omega-1)$, the
the lines of mass $k$ have $y$-intercepts
$(0,\pm k\tau)$.
These equations are all taken mod $\omega$.
\end{lemma}

\begin{theorem}
\label{hier}
Let $B$ be any block.
For each even $k \in [0,p+q]$ there are
$2$ lines in $\cal H$ and $2$ lines
in $\cal V$ which have capacity $k$ and
intersect $B$.  Each such line
carries $k$ light points in $B$.
\end{theorem}

Let $p/q$ be an even rational parameter
and let $\omega=p+q$.  We define the
$x$-diameter of a set to be the
diameter of its projection onto the
$x$-axis. Here is a result repeated from [{\bf S0\/}].

\begin{theorem}
\label{first}
Let $B$ denote the first block.
Then there exists a plaid polygon
in $B$ whose $x$-diameter is at least
$\omega^2/(2q)-1$.  Moreover,
this polygon has bilateral symmetry
with respect to reflection in the
horizontal midline of $B$.
\end{theorem}

\startproof
Let $L$ be the horizontal
line of capacity $2$ and positive sign
which intersects $B$.  Let $z_1=(0,y) \in L$.
By Lemma \ref{anchor}, we
know that $z_1$ is a light point of mass $1$.
Let
$z_2=(\omega^2/2q,y)$.
We compute that $z_2$ is another
light point on $L$.
Since $L$ has capacity $2$, these are
the only two light points on $L$.
The lattice polygon which crosses the
unit horizontal segment containing $z_1$
must also cross the unit horizontal segment
containing $z_2$ because it has to
intersect $L \cap B$ twice.
This gives the lower bound on the $x$-diameter.

Let
$\Gamma'$ denote the reflection
of $\Gamma$ in the horizontal midline of $B$.
We want to show that $\Gamma'=\Gamma$.
Let $V_1$ and $V_2$ denote the two vertical lines
of $B$ having capacity $2$.  These lines are 
symmetrically placed with respect to the
vertical midline of $B$.  Hence, one of the
two lines, say $V_1$, lies less than $\omega /2$
units away from the $y$-axis.
Since 
$\omega /2<\omega^2/(2q),$
the point $z_2$ is separated from the $y$-axis by $V_1$.
Hence both $\Gamma$ and $\Gamma'$ intersect
$V_1$.  Since there can be at most $1$
plaid polygon which intersects $V_1 \cap B$,
we must have $\Gamma=\Gamma'$.
\endproof

\begin{center}
\resizebox{!}{2.6in}{\includegraphics{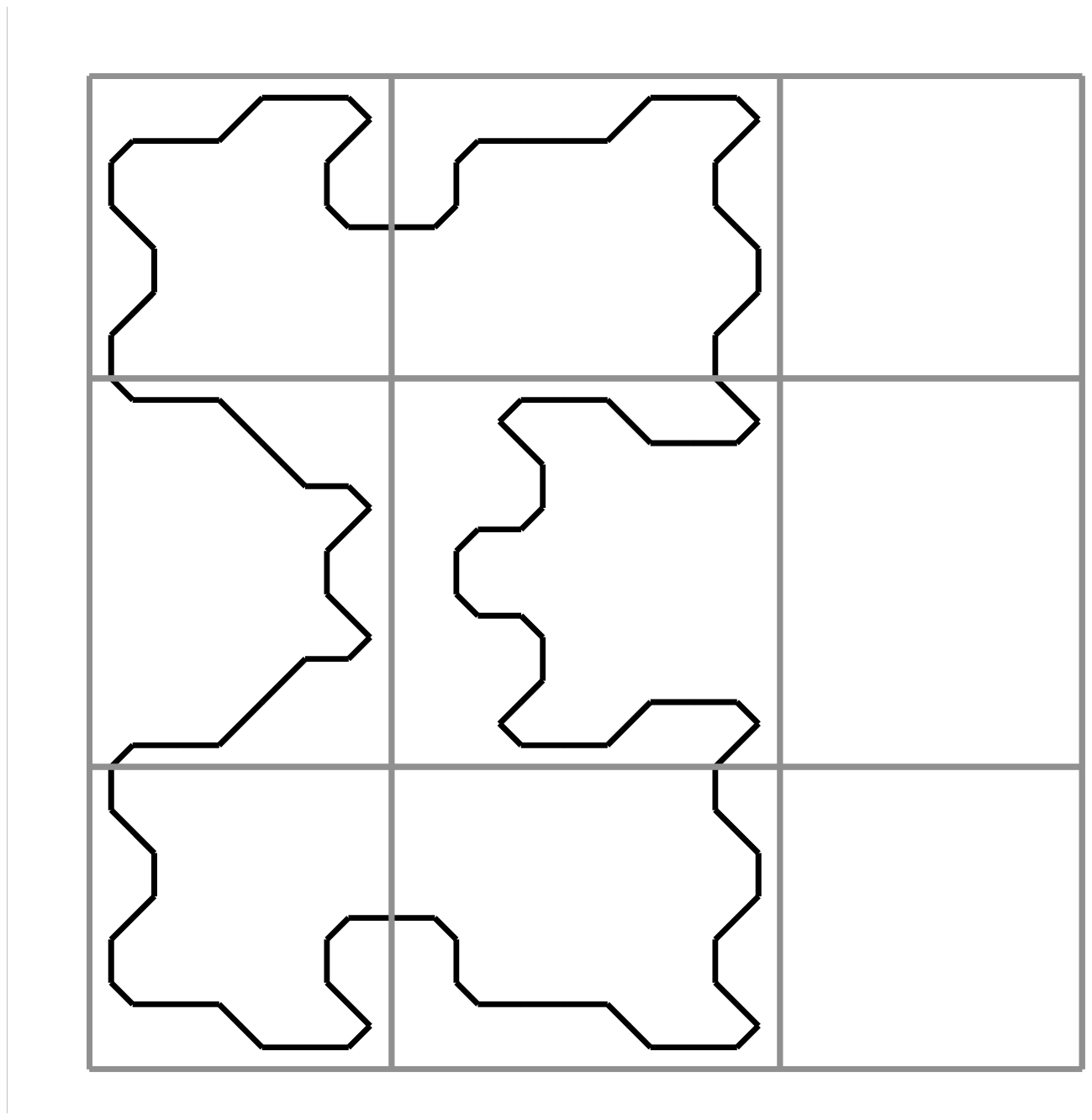}}
\resizebox{!}{2.6in}{\includegraphics{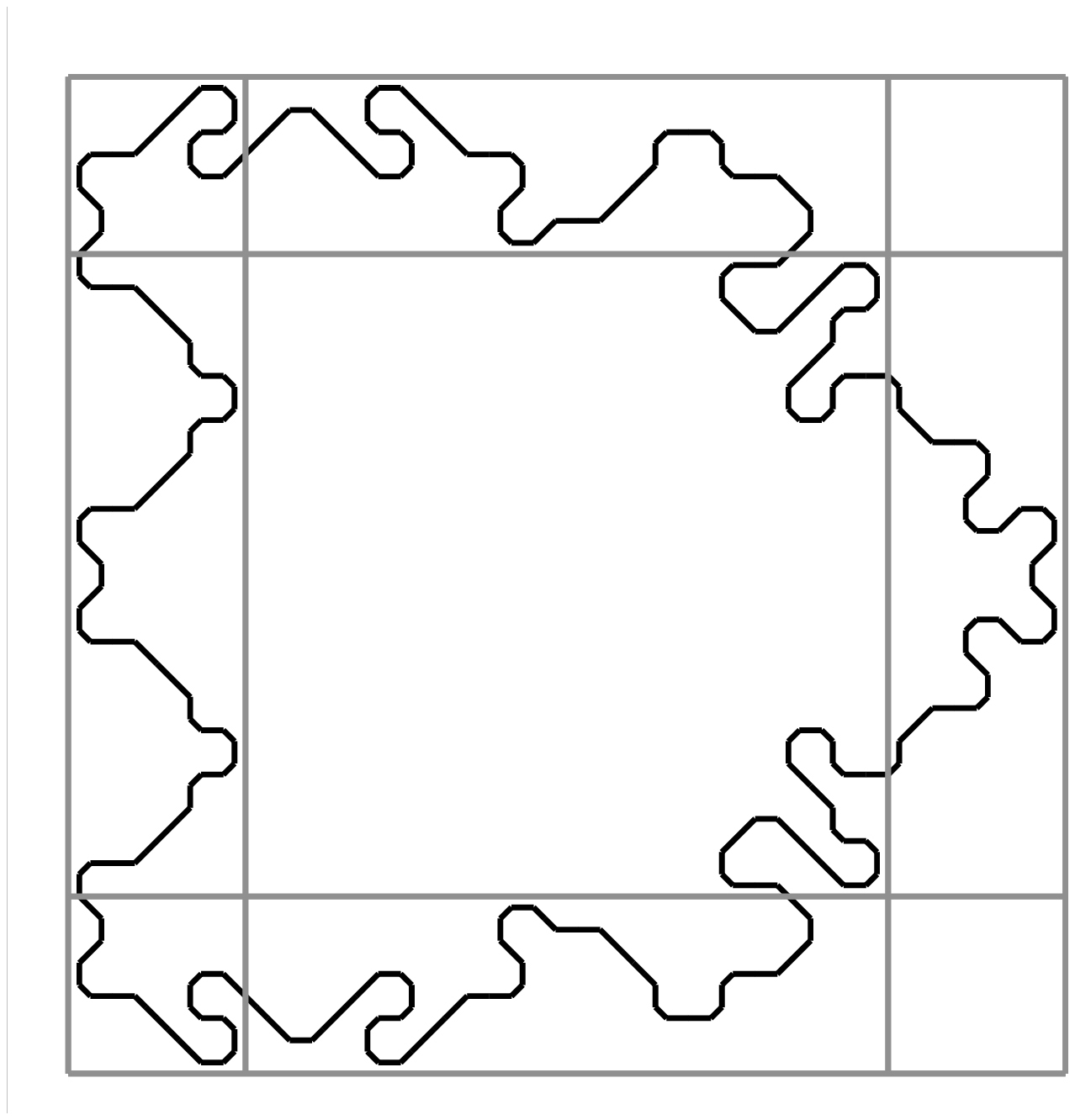}}
\newline
{\bf Figure 2.1:\/} The polygon $\Gamma$ for parameters
$5/18$ and $14/31$.
\end{center}

We call $\Gamma$ {\it the big polygon\/}. Figure 2.1
shows a picture for the $2$ somewhat
arbitrarily chosen parameters.
The lines of capacity $2$ are also shown in
the figure. 

We define
\begin{equation}
\langle\Gamma\rangle=\Gamma \cap (\{1/2\} \times \Z).
\end{equation}
We call points in $\langle \Gamma \rangle$ the {\it anchor points\/}
of $\Gamma$.  We know from our proof that
$\Gamma$ contains at least $2$ anchor points,
namely 
\begin{equation}
\label{special}
z_1+(1/2,0), \hskip 30 pt
z_2+(1/2,0).
\end{equation}

The idea behind our proof of the
Unbounded Orbits Theorem is to take a
geometric limit of a sequence of these
polygons corresponding to a sequence
$\{p_n/q_n\}$ of even rationals
converging to some irrational $A \in (0,1)$. 
As we will explain in \S \ref{fail},
this argument must be done carefully
in order to work.  The key step in our
proof is understanding the set 
$\langle \Gamma \rangle$ well.
We will understand this set in a
recursive way, using the Matching
Criterion to understand what happens
for complicated rational parameters
in terms of simpler ones which are
close in the Diophantine sense.

\newpage

\section{The PET Description}

The material here is a compressed account of
[{\bf S0\/},\S 3] and 
[{\bf S0\/},\S 8].  The main difference
is that the objects we called
$\widehat \Xi$ and $\widehat X$ will
here be called $\Xi$ and $X$.  This
will simplify the notation.

\subsection{The Classifying Space}

We fix $P \in [0,1]$ and let
$\Lambda_P$ denote the lattice generated by the
vectors
\begin{equation}
(0,4,2P,2P), \hskip 30 pt (0,0,2,0), \hskip 30 pt (0,0,0,2).
\end{equation}
This lattice acts on $\R^4$, but
the action on the first coordinate is trivial.

We take the quotient
\begin{equation}
 X_P=(\{P\} \times \R^3)/\Lambda_P.
\end{equation}
The cube
$\{P\} \times [-2,2] \times [-1,1]^2$
serves as a fundamental domain for the action of
$\Lambda_P$ on $\{P\} \times \R^3$. However, the
boundary identifications depend on $P$.
The {\it total space\/} is the
union
\begin{equation}
 X=\bigcup_{P \in [0,1]}  X_P
\end{equation}
This space is a flat affine manifold:  It
is realized as a quotient of the form
\begin{equation}
([0,1] \times \R^3)/\Lambda,
\end{equation}
where $\Lambda$ is the abelian group of
affine transformations generated by
\begin{enumerate}
\item
$T_1(x_0,x_1,x_2,x_3)=(x_0,x_1+4,x_2+4x_0,x_3+4x_0)$.
\item
$T_2(x_0,x_1,x_2,x_3)=(x_0,x_1,x_2+2,x_3)$.
\item
$T_3(x_0,x_1,x_2,x_3)=(x_0,x_1,x_2,x_3+2)$.
\end{enumerate}
The set
$[0,1] \times [-2,2] \times [-1,1]^2$ serves as a
fundamental domain for $ X$.

\subsection{The Partition}

We have a partition of $ X_P$ into
$13=12+1$ regions.  One of the regions is
labeled by an ``empty symbol'' and
the other $12$ regions are labeled by
ordered pairs of elements in $\{N,S,E,W\}$.
For each tile label $\xi$, let $X_P(\xi)$ denote
the piece of the partition labeled by $\xi$.

We describe the partition exactly in
[{\bf S2\/},\S 3].  Here we will just
describe the two relevant features.
First, for each of the $13$ choices
of $\xi$, the {\it total piece\/}
\begin{equation}
X(\xi)=\bigcup_{P \in [0,1]} X_P(\xi)
\end{equation}
is a union of finitely many convex polytopes
whose vertices have integer coordinates.

To describe the second feature, we introduce
coordinates $(P,T,U_1,U_2)$ on $\R^4$.
We have a fibration which maps
$(P,T,U_1,U_2)$ to $(P,T)$.  The fibers
are called the $(U_1,U_2)$ fibers.  The second
feature is that the partition intersects
each $(U_1,U_2)$ fiber generically in a
$4 \times 4$ grid of rectangles.  
The rectangles are parallel to the coordinate
axes.  Figure 3.1 shows what we have in mind.
The squares labeled N,E,W,S are really labeled
by the empty symbol in our scheme.  However,
we put these letters in them do indicate the
scheme for filling in the rest of the labels.

\begin{center}
\resizebox{!}{3in}{\includegraphics{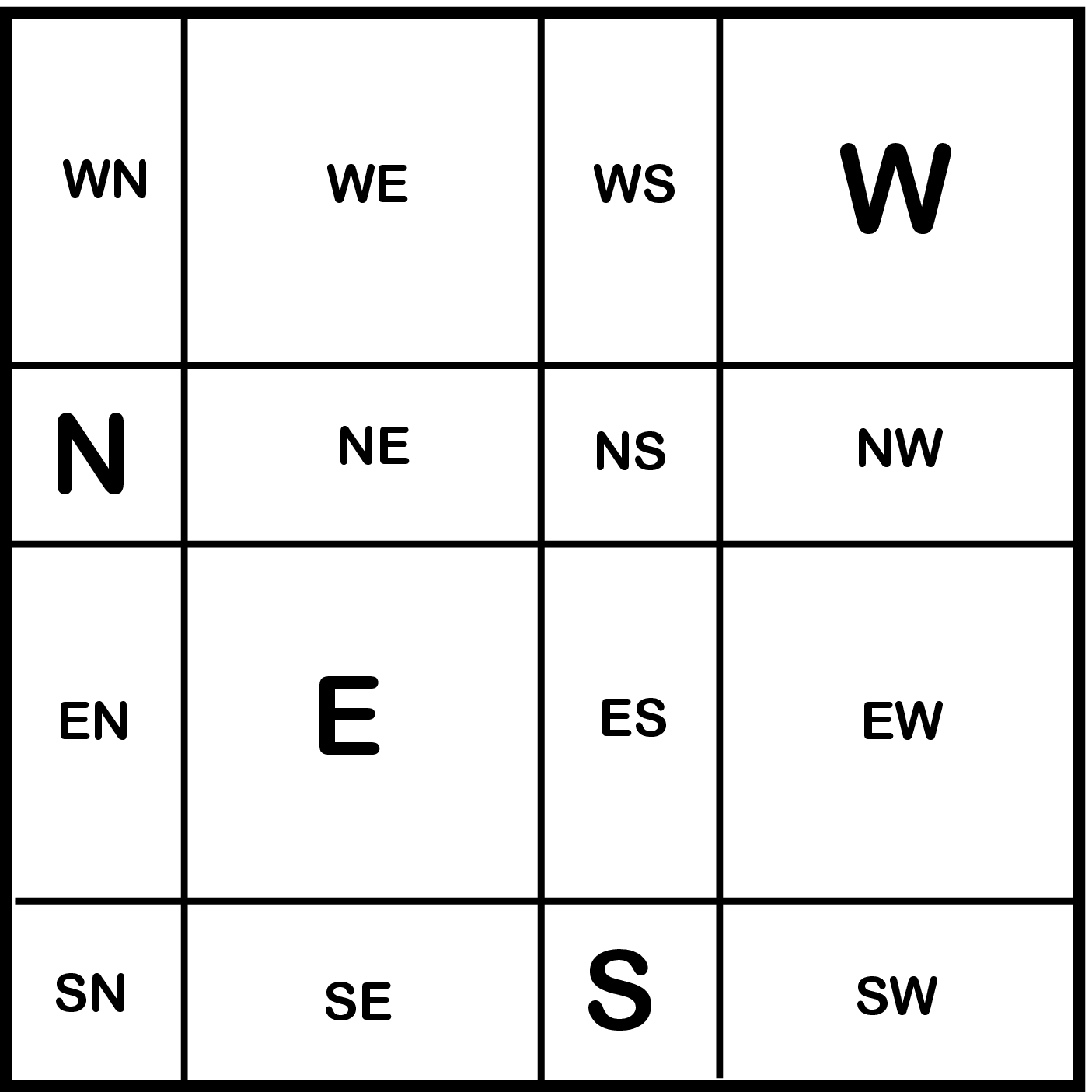}}
\newline
{\bf Figure 3.1:\/} The intersection of the
partition with a fiber.
\end{center}

\subsection{The Classifying Map}

For each $P \in (0,1)$ we have a map
\begin{equation}
\Xi_P: \R^2 \to X_P
\end{equation}
defined by
\begin{equation}
\label{map}
\Xi(x,y) = (P,2Px+2y,2Px,2Px+2Py)
\hskip 10 pt {\rm mod\/}\hskip 5 pt \Lambda_P.
\end{equation}
For each even rational parameter $p/q$ we
set $P=2p/(p+q)$, as above, and
restrict $\Xi_P$ to the set
$\cal C$ of centers of square tiles.
In [{\bf S0\/}] we proved that
$\Xi_P({\cal C\/})$ never hits a
wall of a partition in the rational
case. 
 
For fixed $p/q$, we define a tiling of
$\R^2$ as follows.  For each $c \in \cal C$
we see which partition piece contains the
point $\Xi$, and then we draw the corresponding
tile in the unit square centered at $c$.
If $\Xi_P(c)$ lands in the tile labeled
by the empty symbol, then we draw nothing in
the tile.

The PET Equivalence Theorem in
[{\bf S0\/}] says that
the tiling we get agrees with the tiling 
from the plaid model, and moreover
the orientations on the tiles are cpherent.
That is, for each edge involved in two connectors,
one of the connectors points into the edge and
the other points out.  That is, each of the
plaid polygons has a coherent orientation.

The next result says something about the
geometry of the clasifying map.  The result
holds for any half integer, but we only care
about the case of $1/2$.

\begin{lemma}
\label{geodesic0}
If $Y=\{1/2\} \times \Z$ then
$\Xi_P(Y)$ is contained in a single geodesic
$\gamma_P$ in 
the fiber of $\Xi_P$ over the point $(P,P)$.
Moreover, $\gamma_P$ has slope $-1$ when developed
into the plane so that the $(U_1,U_2)$ coordinate
axes are identified with the coordinate
axes of the plane.
\end{lemma}

\startproof
  To see this, we compute
$$
\Xi_P(1/2,m+1/2)=(P,P+2m,P,P+2Pm) \hskip 10 pt
{\rm mod\/} \hskip 10 pt \Lambda_P.
$$
Subtracting off the vector $(2m,Pm,Pm)$, we get
\begin{equation}
\label{geodesic}
\Xi_P(1/2,m+1/2)=(P,P,P,P)+(0,0,-mP,mP)  \hskip 10 pt
{\rm mod\/} \hskip 10 pt \Lambda_P.
\end{equation}
As $m$ varies, the image remains on the geodesic
of slope $-1$ through the point $(P,P,P,P)$.
At the same time the walls of the partition
intersect the fiber in horizontal and vertical lines.
\endproof

\subsection{The PET Interpretation}

Now we see how to interpret the
space $ X_P$ as a PET

For each parameter $p/q$, we set $P=2p/(p+q)$, 
and our space 
is $ X_P$.
The partition is given by
\begin{equation}
 X_{S \downarrow} \cup
 X_{W \leftarrow} \cup
 X_{N \uparrow} \cup
 X_{E \rightarrow}.
\cup  X_{\Box}.
\end{equation}
The remaining set $ X_{\Box}$ is just the complement.
These sets have an obvious meaning.  
For instance, $ X_{E,\rightarrow}$ is the union of
regions which assign tiles which point into their east edges.
And so on.
The second partition of $ X$ is obtained by
reversing all the arrows.  So, the first partition
is obtained by grouping together all the regions which
assign tiles which point into a given edge, and
the second partition is obtained by grouping together
all the regions which assign tiles which point
out of a given edge.

We have $5$ {\it curve following maps\/}
$\Upsilon_{\Box}:  X_{\Box} \to  X_{\Box}$ is
the identity, and then
$\Upsilon_S:  X_{S \downarrow} \to  X_{N\downarrow}$.
And so on.
These maps are all translations on their domain.
For instance, $\Upsilon_P$ has the following property:
\begin{equation}
\Upsilon_S \circ \Xi(x,y)=\Xi(x,y-1). 
\end{equation} 
From Equation \ref{map}, we get
\begin{equation}
\Upsilon_S(P, T,U_1,U_2)=(P, T-2,U_1,U_2-2P) \hskip 10 pt
{\rm mod\/} \hskip 10 pt \Lambda'.
\end{equation}

There is a natural map on $\cal C$, the set of centers
of the square tiles.  We simply follow the directed edge
of the tile centered at $c$
and arrive at the next tile center.  By construction,
 $ \Xi_P$ conjugates this map on
$\cal C$ to our polyhedron exchange transformation.
By construction $ \Xi_P$
sets up a dynamics-respecting bijection 
between certain orbits in 
$ X_P$ and the plaid polygons
contained in the fundamental domain, with
respect to the parameter $p/q$.

Now we explain what we mean by vector dynamics.
We assign
$(0,-1)$ to the region $ X_{S\downarrow}$, 
and $(-1,0)$ to $ X_{W\leftarrow}$, etc.
When we follow the orbit of
$\Xi(c)$, we get record the list of vectors labeling
the regions successively visited by the point. This gives
vectors $v_1,v_2,v_3$.  The vectors
$c$, $c+v_1$, $c+v_1+v_2$, etc. are the vertices
of the plaid polygon containing $c$.
This is what we mean by saying that the plaid polygons
describe the vector dynamics of certain orbits of
the PET.

\subsection{The Irrational Case}
\label{irr}

Let $\cal C$ denote the set of
centers of the square tiles.
The classifying pair $(\Xi_P, X_P)$ makes
sense even when $P$ is irrational.
However, $ \Xi_P$ might map
points of $\cal C$ into the boundary
of the partition.
For instance
\begin{equation}
\Xi_P(1/2,1/2)=(P,P+1,P,2P)
\end{equation}
and this point lies in the wall of
the partition.

To remedy this situation, we introduce an
{\it offset\/}, namely a vector
$V \in \R^3$, and we define
\begin{equation}
\Xi_{P,V}=\Xi_P+(0,V).
\end{equation}
(The first coordinate, namely $P$, does not change.)
Since $\Lambda_P$ acts on $\R^3$ with
compact quotient, we can take $V$ to lie in a
compact subset of $\R^3$ and still we will
achieve every possible map of this form.
We call $V$ a {\it good offset\/} if
$\Xi_{P,V}(\cal C)$ is disjoint from
the walls of the partition.
In this case, the PET dynamics
permutes the points of 
$\Xi_{P,V}(\cal C)$ and every orbit of one
of these points is well defined.

Here we give a criterion for $V$ to be a good offset.
Define
\begin{equation}
\Q[P]=\{r_1+r_2P|\ r_1,r_2 \in \Q\}.
\end{equation}

\begin{lemma}
\label{good0}
Suppose $V=(V_1,V_2,V_3)$ is such that
$$V_1 \in \Q[P], \hskip 30 pt
V_1,V_2 \not \in \Q[P].$$ Then $V$
is a good offset.  
\end{lemma}

\startproof
Here is an immediate consequence the two main
properties of our partition.  Suppose
$F$ is a $(U_1,U_2)$ fiber over a point
$(P,T) \in \Q[P] \times \Q[P]$.  Then
the walls of the partition intersect $F$ in
rectangles bounded by lines defined
over $\Q[P]$.  That is, the lines have
the form $x=x_0$ and $y=y_0$ for various
choices of $x_0,y_0 \in \Q[P]$.  So,
if we have a point $(P,T,U_1,U_2)$ with
$T \in \Q[P]$ and $U_1,U_2 \not \in \Q[P]$,
then the point lies in the interior of
some partition piece.  
Given the formula in Equation \ref{map}, the
map $ \Xi_{P,V}$ has this property.
\endproof

\subsection{Geometric Limits}

Suppose now that $\{p_k/q_k\}$ is a sequence
of even rational parameters converging to
some irrational parameter $A$. 
Let 
\begin{equation}
P_k=\frac{2p_k}{p_k+q_k},
\hskip 30 pt P=\frac{2A}{1+A}.
\end{equation}
By definition $P_k \to P$.  Likewise,
the classifying maps
$ \Xi_{P_k}$ converge, uniformly
on compact sets, to $ \Xi_P$.
As we have already mentioned, the map
$ \Xi_P$ is not really a very
good map.  Here we explain a more flexible
kind of limit we can take.

For ease of notation, we set 
$ \Xi_{P_k}= \Xi_k$.
Let $Z=\{z_k\}$ be a sequence of points
in $\Z^2$.  We define
\begin{equation}
 \Xi_{k}^Z(c)= \Xi_{k}(c+z_k).
\end{equation}
In other words, we compose by a translation, so
that the action of $ \Xi_{k}^Z$ around
the origin looks like the action of
$ \Xi_{k}$ near $z_k$.

\begin{lemma}
The sequence $\{ \Xi_k^Z\}$ converges on
a subsequence to
$ \Xi_{P,V}$ for some offset vector $V$.
\end{lemma}

\startproof
Since $ \Xi$ is a locally affine map
for each parameter, there is some vector
$V_k$ so that
$$ \Xi_k^Z =  \Xi_k + V_k.$$
We can always take $V_k$ to lie within a
compact set, 
namely $$\{0\} \times [-2,2] \times [-1,1]^2,$$
which is just a translated copy of the fundamental
domain for the lattice $\Lambda_{P_k}$. (It is
translated so that the first coordinate is $0$.)
But now we can take a subsequence so that the
sequence $\{V_k\}$ converges to some vector $V$.
By construction, 
 $\{ \Xi_k^Z\}$
converges to $\Xi_{P,V}$.
\endproof

The vector $V$ depends somehow on the sequence $Z$.
We call $Z$ a {\it good sequence\/} if the
vector $V$ is a good offset.

\subsection{A Failed Attempt}
\label{fail}

To explain the subtlety in the
proof of the Unbounded Orbits Theorem,
we make a failed attempt at proving it.

We choose 
some sequence $\{p_k/q_k\}$ which converges
to $A$ and let $\{\Gamma_k\}$ be the
corresponding sequence of big polygons.
Let $z_k$ be a point of $\cal C$ uniformly
close to one of the two points
of $\langle \Gamma_k \rangle$ that
we know about already, namely
those in Equation \ref{special}.
Let $Z=\{z_k\}$.  

Now, $\Gamma_k-z_k$ intersects a
tile centered at a point $c_k$ which is
uniformly close to the origin.  Hence,
on a subsequence, the polygons
$\{\Gamma_k-z_k\}$ converge to an infinite
polygonal path $\Gamma_{\infty}$.
If $Z$ is a good sequence, then
$\Xi_{P,V}=\lim \Xi_k^Z$
has an infinite orbit
whose vector dynamics produces
$\Gamma_{\infty}$.  We will
discuss this point in more detail 
in \S \ref{endproof}.

One problem with this construction is that
we don't know that
both coordinate projections of
$\Gamma_{\infty}$ have infinite
diameter. 
A more serious problem is that
$Z$ is not a good sequence.
A calculation, which we omit, shows
that $\Xi_P$ maps the points $z_j+(\pm 1/2,\pm 1/2)$
into the walls of the partition.

In making an argument which works, we want to
keep the idea of choosing points
near $\langle \Gamma_k \rangle$, but
we need to know that this set has many
more points before we can make
fundamentally different choices.

\newpage

\section{The Main Argument}

\subsection{Approximating Irrationals}
\label{strong}

Our goal in this section is to define
a sequence $\{p_n/q_n\}$ of even
rationals which converges to an
irrational $A \in (0,1)$ and has
good Diophantine properties.
This sequence is related to, but
usually different from, the sequence
of continued fraction approximants.
We defer the proofs of the technical
lemmas to the next chapter.
\newline
\newline
{\bf Even Predecessors:\/}
Two rationals $a_1/b_1$ and $a_2/b_2$ are
called {\it Farey related\/} if
$|a_1b_2-a_2b_1|=1$.
Let $p'/q'$ and $p/q$ be two even rational
paramters.  We write
$p'/q' \leftarrow p/q$ if $p/q$ and $p'/q'$ are
Farey related and
$\omega'<\omega$.  Here $\omega=p+q$ as usual.
We call $p'/q'$ the {\it even predecessor\/} of $p/q$.
This rational is unique.  
\newline
\newline
{\bf Core Predecessors:\/}
Recall that $\tau$ is the tune of $p/q$,
defined in \S \ref{tune}.
We define $\kappa \geq 0$ to be the integer so that
\begin{equation}
\frac{n}{2n+1} \leq \frac{\tau}{\omega} < \frac{n+1}{2(n+1)+1}.
\end{equation}
We only get equality on the left hand side when
$p/q=1/2n$.
Define 
\begin{equation}
\label{corep}
\widehat p=p-2\kappa p', \hskip 30 pt
\widehat q=q-2\kappa q'.
\end{equation}
Statement 1 of Lemma \ref{omnibus} proves
that $\widehat p/\widehat q$ is an
even rational in $(0,1)$.
We call $\widehat p/\widehat q$ the
{\it core predecessor\/}.
\newline
\newline
\noindent
{\bf Predecessors:\/} Given an even rational
parameter $p/q$, we define the {\it predecessor\/}
$p^*/q^*$
as follows:
\begin{itemize}
\item If $p=1$ then $p^*/q^*=0/1$.
\item If $p \geq 2$ and $\kappa=0$ then $p^*/q^*=p'/q'$,
the even predecessor of $p/q$.
\item If $p \geq 2$ and $\kappa \geq 1$ then $p^*/q^*=\widehat p/\widehat q$,
the core predecessor of $p/q$.
\end{itemize}
We write $p^*/q^* \prec p/q$.
This definition turns out to be very well adapted to
the plaid model.  Lemma \ref{omnibus} collects together
many of the relations between these rationals.
\newline
\newline
{\bf The Predecessor Sequence:\/}
In \S 5 we will prove the following result.
\begin{lemma}
\label{neat}
Let $A \in (0,1)$ be irrational.
Then there exists a sequence
$\{p_k/q_k\}$ such that \begin{itemize}
\item $p_0/q_0=0/1$
\item $p_k/q_k \prec p_{k+1}/q_{k+1}$ for all $k$
\item $A=\lim p_k/q_k$.
\end{itemize}
\end{lemma}
We call $\{p_k/q_k\}$ the {\it predecessor sequence\/}.
This terminology suggests that the predecessor
sequence is unique.  However, we will not
bother to prove this.  We just need existence,
not uniqueness. 
\newline
\newline
{\bf Diophantine Result:\/}
Let $\{p_k/q_k\}$ be the predecessor sequence
converging to $A$. 
We classify a term $p_k/q_k$ in the
predecessor sequence as follows:
\begin{itemize}
\item weak: $\tau_{k+1}<\omega_{k+1}/4$.  Here $\kappa_{k+1}=0$.
\item strong: $\tau_{k+1} \in (\omega_{k+1}/4,\omega_{k+1}/3)$.
Here $\kappa_{k+1}=0$.
\item core: $\tau_{k+1}>\omega_{k+1}/3$.
Here $\kappa_{k+1} \geq 1$.
\end{itemize}

\begin{lemma}
\label{strong0}
The predecessor sequence has infinitely many non-weak
terms.  For each non-weak term $p_k/q_k$, we have
$$\bigg|A-\frac{p_k}{q_k}\bigg|<\frac{48}{q_k^2}.$$
\end{lemma}

\noindent
{\bf The Approximating Sequence:\/}
We define the approximating sequence to be
the set of terms $p_k/q_k$ in the predecessor sequence
such that either
\begin{itemize}
\item $p_k/q_k$ is core.
\item $p_k/q_k$ is strong and
$p_{k-1}/q_{k-1}$ is not core.
\end{itemize}
If there are infinitely many core terms in
the predecessor sequence, then the approximating sequence
contains all of these.  If there are only finitely
many core terms in the predecessor sequence, then
there are infinitely many strong terms, and the
approximating sequence contains all but finitely many
of these.  So, in all cases, the approximating sequence
is an infinite sequence.

\subsection{Arc Copying}
\label{mb}

Let $R_{p/q}$ denote the rectangle bounded by
the bottom, left, and top of the first block,
and by whichever vertical line of capacity
at most $4$ is closest to the left edge of
the first block.
Let $\gamma_{p/q}$ denote the subset of $\Gamma_{p/q}$
contained in the box $R_{p/q}$.  When the dependence
on the parameter is implied we will suppress it
from our notation.

Figure 4.1 shows the big polygons
associated to two different parameters.
Notice that the polygon $\Gamma_{12/29}$
copies some of $\Gamma_{5/12}$.  The
boxes $B_{5/12}$ and $B_{12/29}$ are
the first columns (i.e. the union of the
leftmost $4$ sub-rectangles) of the tic-tac-toe
grids shown in each of the two pictures.

\begin{center}
\resizebox{!}{1.3in}{\includegraphics{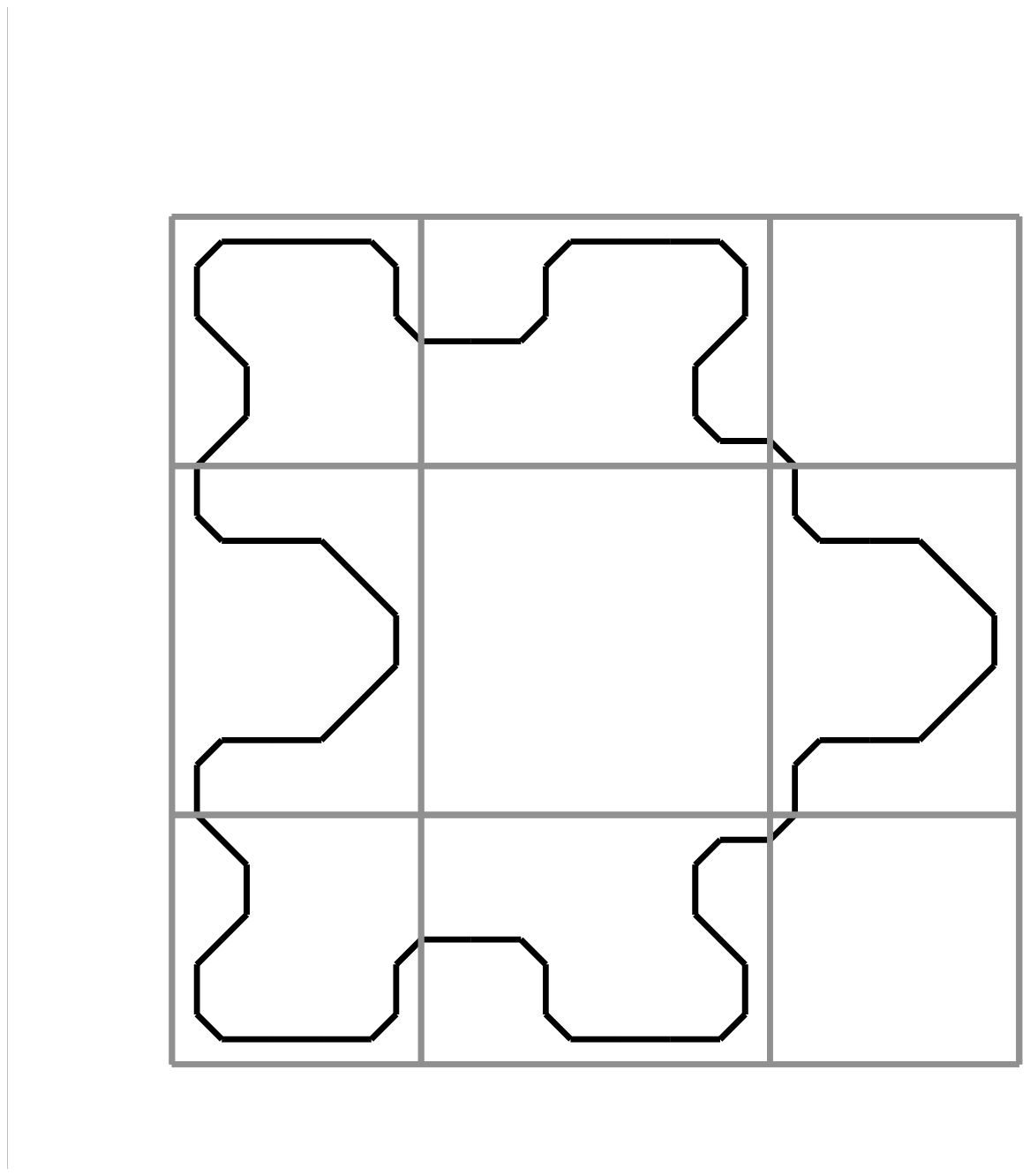}}
\resizebox{!}{2.6in}{\includegraphics{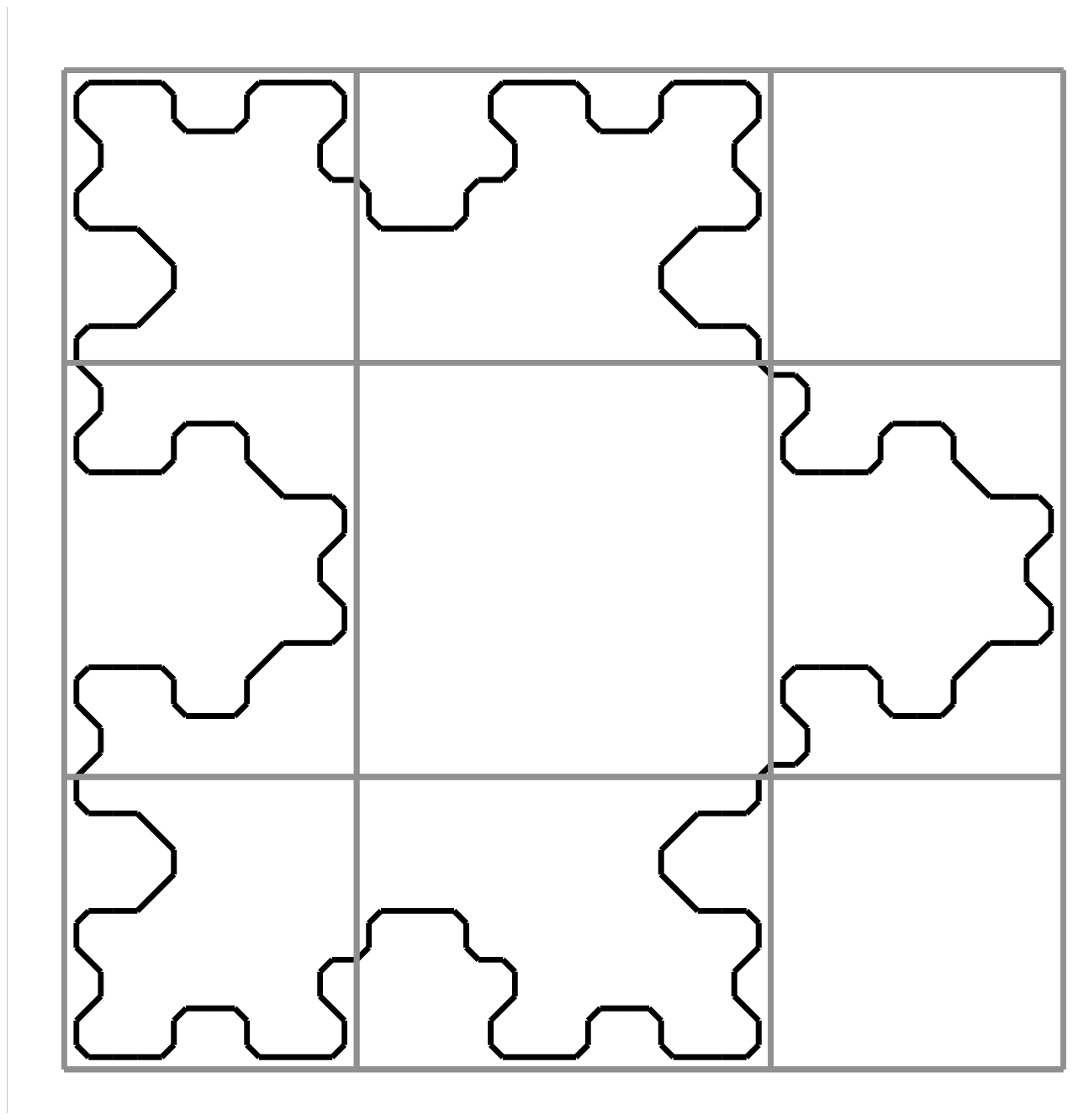}}
\newline
{\bf Figure 4.1:\/} Arc copying for $5/12$ and $12/29$.
\end{center}

Let $TH$ and $BH$ denote the top and bottom horizontal
lines of capacity $2$ 
with respect to some rational parameter.
In the statement of the results below, it will be clear
which parameters these lines depend on.

\begin{lemma}[Box]
For any even rational parameter $p/q$, the set
$\gamma_{p/q}$ is an arc whose endpoints
lie on the right edge of $R_{p/q}$.
\end{lemma}

\begin{theorem}[Copy]
Let $p_0/q_0, p_1/q_1$ be two
successive terms in the approximating sequence.
Then there is some vertical translation $\Upsilon$ such that
that $\Upsilon(R_0)$ is contained below the
horizontal midline of $R_1$, and
$T(\gamma_0) \subset \gamma_1$.
Moreover, either $\Upsilon(BH_0)=BH_1$ or
$\Upsilon(TH_0)=BH_1$.
\end{theorem}

\subsection{Marked Boxes}

Now we will
somewhat abstract the main features of the
pair $(R_{p/q},\gamma_{p/q})$, even
though approximatingly these are the sets
we will be considering.

Say that a {\it marked box\/} is a pair
$\beta=(R,\gamma)$
\begin{itemize}
\item $R$ is a rectangle with sides parallel to the
coordinate axes, whose left edge is contained 
in the $y$ axis. 
\item $\gamma$ is a polygonal path which
has both endpoints on the right edge of $R$ and
contains two points of the form
$(1/2,y)$, one strictly in the lower half of $R$
and one strictly in the upper half.
\end{itemize}

Given two marked boxes
$\beta_j=(R_j,\gamma_j)$, we write
$\beta_1 \prec \beta_2$ if 
$R_1 \subset R_2$ and $\gamma_1 \subset \gamma_2$.
We insist that the side lengths of $R_2$ are
at least one unit longer than the corresponding
side lengths of $R_1$.  This is the
kind of copying produced by the
Copy Theorem.

Now we explore the consequences of
iterating the result of the
Copy Theorem.
Note that
the Copy Theorem really says that the arc $\gamma$
contains two translated copies of $\widehat \gamma$,
one above the horizontal midline of $R$ and one below.
This follows from the bilateral symmetry
of $\gamma$ and of $\widehat \gamma$.
So, when we iterate the Copy Theorem,
we should expect the emergence of
a tree-like structure.

Let $T$ be a directed binary tree with no
forward infinite path.  If $T$ is finite,
then $T$ has one initial vertices
$2^k-1$ vertices in total.
If $T$ is infinite, then $T$ has no initial
vertex. Given any vertex $v$ of $T$,
we define $\wedge v$ to be the subtree of
$T$ whose initial node is $v$.  This is
always a finite binary tree.

We say that a 
collection of marked boxes {\it realizes\/}
$T$ if there is map $\Phi$ from the set of
vertices of $T$ to the set of marked
boxes, such that
\begin{itemize}
\item If $v \leftarrow w$ then $\Phi(v) \prec \Phi(w)$.
\item If $v_j \leftarrow w$ for $j=1,2$, then
$\Phi(v_1)$ and $\Phi(v_2)$ are disjoint and
the restriction of $\Phi$ to $\wedge v_1$ equals
the restriction of $\Phi$ to $\wedge v_2$ up
to a vertical translation $\tau_w$.  We call
$\tau_w$ the {\it translation associated to\/} $w$.
\end{itemize}

Suppose we have a finite
string $p_1/q_1,...,p_n/q_n$ of parameters in
the approximating sequence.  The Box Lemma says that all
the pairs $(R_k,\gamma_k)$ are marked boxes.
Iterated application of
the Copy Theorem,
and bilateral symmetry,
gives us a realization of the directed binary
tree with $2^n-1$ vertices.

Figure 4.2 shows a schematic picture of the
case when $n=3$.

\begin{center}
\resizebox{!}{3in}{\includegraphics{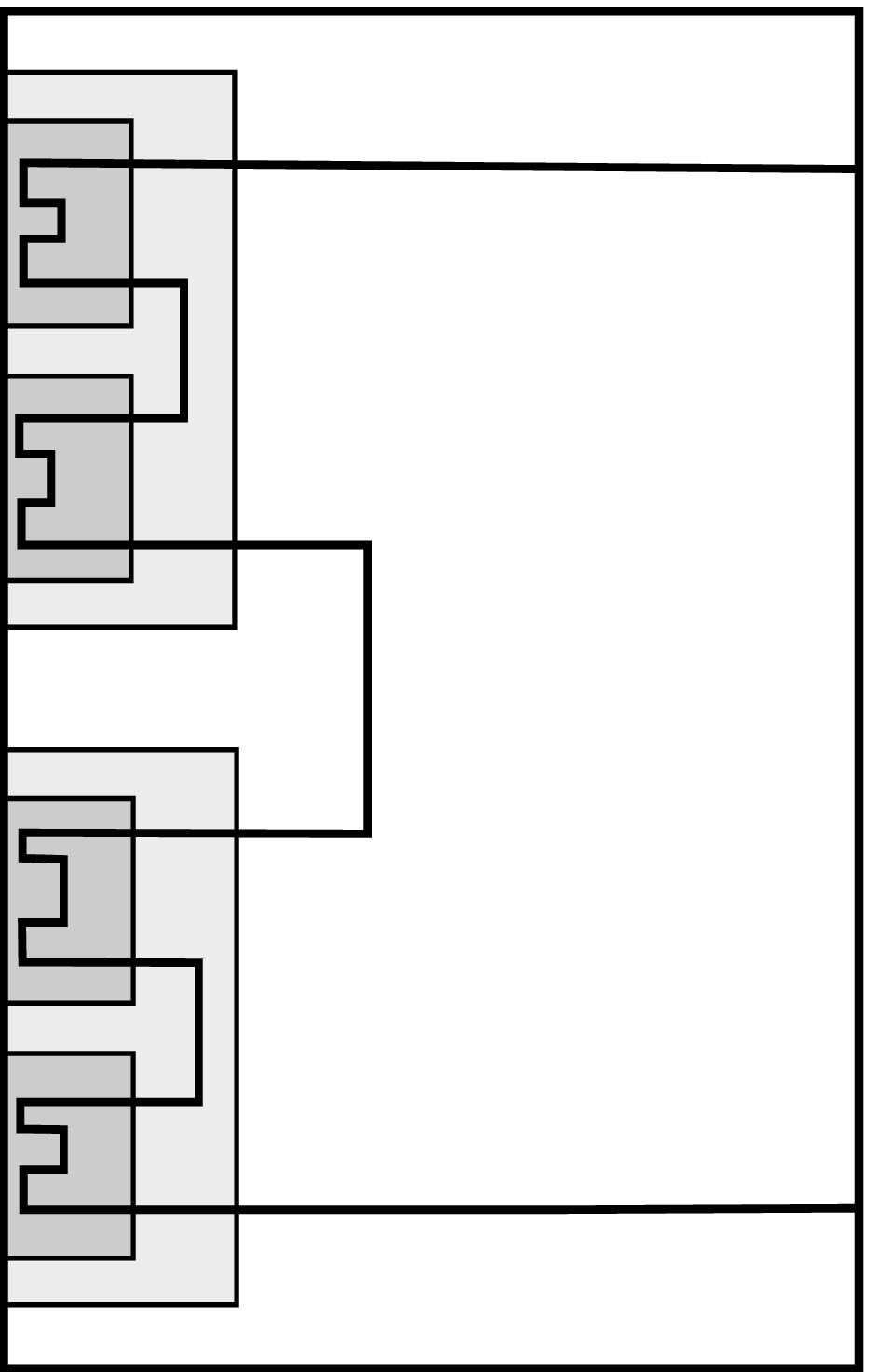}}
\newline
{\bf Figure 4.2:\/} A pattern of boxes realizing a tree.
\end{center}

Figure 4.3 shows a particular planar
embedding of an infinite binary directed
tree.  The bottom path lies in the $x$-axis.
The terminal nodes of the tree lie on the $y$-axis.

\begin{center}
\resizebox{!}{1.8in}{\includegraphics{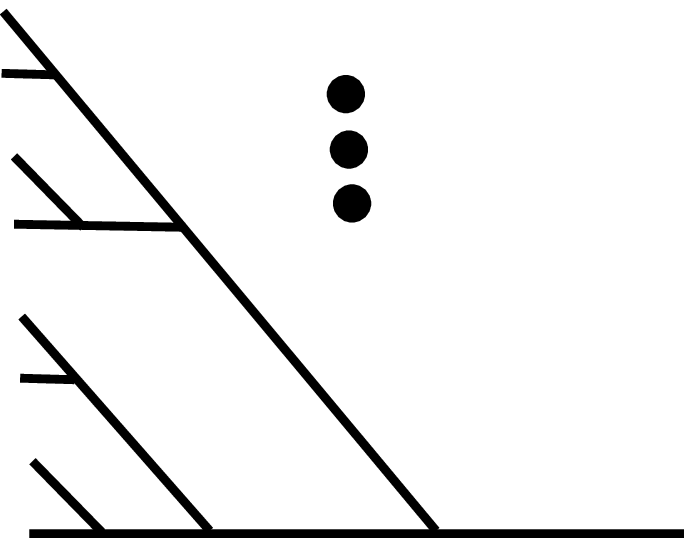}}
\newline
{\bf Figure 4.3:\/} A special infinite binary tree.
\end{center}

Let $\Phi$ be a realization of the tree $T$
suggested by Figure 4.3.
We call $\Phi$ {\it normalized\/} if
$\Phi$ maps any vertex on the bottom path of
$T$ to a box whose curve contains $(1/2,0)$,
and if the boxes obtained in this way are
beneath any other boxes associated to $\Phi$.
This is to say that all the translations
$\tau_w$ considered in the definition of $\Phi$
are positive translations along the $x$-axis.

\subsection{The Realization Lemma}

Let $T$ be an infinite directed binary tree,
as above.  

\begin{lemma}[Realization]
Let $A \in (0,1)$ be irrational.
Let $\{p_n/q_n\}$ be the approximating sequence
associated to $A$.
Then there exists a normalized
tree realization $\Phi$ of the
infinite binary tree $T$ with
the following property.
For each $n$ there is a vertex
$v_n$ of $T$ such that
$\Phi(v_n)$ is translation equivalent
to the marked box $\beta_{p_n/q_n}$.
\end{lemma}

\startproof
Looking at the finite sequence
$p_1/q_1,...,p_n/q_n$ of parameters
in the approximating sequence,
we get a realization $\Phi_n$ of the tree
$T_n$ having $2^n-1$ vertices.  We can
adjust $\Phi_n$ by a vertical translation so that
the bottom distinguished point of $\Phi_n$ is
$(1/2,0)$.  By construction, the bottom half
of $\Phi_{n+1}$ coincides with $\Phi_n$.
This means that we can take a geometric limit
and get a normalized realization of the
infinite tree.  The vertex sequence
$\{v_n\}$ is just the sequence of vertices
along the $x$-axis in the tree shown in Figure 4.3.
\endproof

Let $d_k$ be the translation length
of the translation $\tau_{v_k}$ associated
to the realization $\Phi$.
Recall that $[t]_2$ denotes the value of
$t$ mod $2\Z$ which lies in $[-1,1)$.
Let $P=2A/(1+A)$.
\begin{lemma}
\label{adapt}
$\lim_{k \to \infty} [Pd_k]_2=0$.
\end{lemma}

\startproof
Let $t_k=t(p_k/q_k)$.  Let
$\eta_k=\omega_k-2t_k$.
Geometrically, $\eta_k$ is the distance between
the two horizontal lines of capacity $2$ associated
to $p_k/q_k$.  It follows from the
last Statement of the Copy Theorem that
$d_k=\eta_k \pm \eta_{k-1}$.
Thus, it suffices to prove that
$\lim_{k \to \infty} [P\eta_k]_2=0$.

Let $P_k=2p_k/(p_k+q_k)$.  The function
$y \to [P_ky]_2$ assigns value of
$\pm 1/\omega_k$ to points on the two
horizontal lines of capacity $2$.  Hence
\begin{equation}
\label{dio2}
\lim_{k \to \infty} [P_k\eta_k]_2 \to 0.
\end{equation}

\begin{equation}
\label{dio3}
|(P-P_k)\eta_k| \leq |P-P_k|\omega_k<96/\omega_k.
\end{equation}
The last inequality comes from Lemma \ref{strong0}
and the fact that the map $A \to P=2A/(A+1)$
is $2$-lipschitz.
This lemma now follows from the triangle
inequality and Equations \ref{dio2} and
\ref{dio3}.
\endproof

\subsection{Taking Many Limits}

Let $\Phi$ be the normalized tree
realization from The Realization Lemma.
Let $\{d_k\}$ be the sequence of
translations associated to $\Phi$, as
in Lemma \ref{adapt}.
Let $B=\{\epsilon_k\}$ be any infinite
binary sequence.  Define
\begin{equation}
\label{sigma}
\Phi_k=\Phi-(0,y_k), \hskip 30 pt
y_k=\sum_{i=1}^k \epsilon_i d_i.
\end{equation}
Note that $(0,1/2)$ remains a vertex
of some of the curves associated to
$\Phi_k$. 

\begin{lemma}
\label{geolim}
For any positive $R$, the images of
$\Phi_k$ and $\Phi_{k+1}$ agree inside the
ball of radius $R$ about $0$ for all
sufficiently large $k$.
\end{lemma}

\startproof
The realizations $\Phi_k$ and $\Phi_{k+1}$
either agree, or there are three
vertices $v_1,v_2,w$ such that
\begin{itemize}
\item $v_1,v_2 \leftarrow w$ and $v_1$ lies below $v_2$.
\item Both $\Phi_k(w)$ and $\Phi_k(v_1)$ contain $(0,1/2)$.
\item $\Phi_{k+1}=-\tau_w \circ \Phi_k$.
\end{itemize}
Moreover, the vertex $w$ is combinatorially $k$
steps from a terminal vertex.  We have
$$\Phi_{k+1}(\wedge v_2)=-\tau_w \Phi_k(\wedge v_2)=
\Phi_k(\wedge v_1).$$
Once $k$ is sufficiently large, 
$\Phi_k(\wedge v_1)$ and $\Phi_{k+1}(\wedge v_2)$
account for all the boxes associated to
these realizations which come within $R$ of the origin.
\endproof

Lemma \ref{geolim} allows us to define the limit
$\Phi_{\infty}=\lim_k \Phi_k$.
By construction, $\Phi_{\infty}$ is another realization
of an infinite tree.  The curve associated
to $\Gamma_{\infty}$ -- namely, the union of all
the curves associated to all the boxes -- has 
infinite diameter projections
onto both coordinate axes.

Let $y_k$ be as in
Equation \ref{sigma}.  Let
$\Xi_{P,V}$ be the a limit of the
translated classifying maps
\begin{equation}
\Xi_P-(0,y_k).
\end{equation}
This definition depends on the binary
sequence $B=\{\epsilon_k\}$.
We call $B$ a {\it good sequence\/} if $V$ is a
good offset.  As we have already mentioned,
one way to guarantee this condition is to
arrange that $\Xi_{P,V}(1/2,1/2)$ has
coordinates of the form $(P,T,U_1,U_2)$ with
$T \in \Q[P]$ and $U_1,U_2 \not \in \Q[P]$.

\begin{lemma}
There is a
good binary sequence relative to $\Phi$.
\end{lemma}

\startproof
Let $Y$ and $\Xi_P$ and $\gamma_P$ be
the sets from Lemma \ref{geodesic0}.
Since $\gamma_P$ has slope $-1$ when
developed out into the plane, the
geodesic $\gamma_P$
only intersects the set
\begin{equation}
(\Q[P] \times \R) \cup (\R \times \Q[P])
\end{equation}
in a countable collection of points.

Let $Y_D \subset Y$ denote the set of all points of the form
\begin{equation}
\bigg(\frac{1}{2},\sum_{i=1}^k \epsilon_i d_i\bigg).
\end{equation}
Here the sum ranges over all finite binary
sequences, and all $k \in \N$.
To prove our result, it suffices to prove that
$\Xi_P(Y_D)$ has uncountably many accumulation
points contained inside $\gamma_P$.
Since $\gamma_P$ is dense in the
fiber containing it, we have to
take care that the
accumulation points we find all
belong to $\gamma_P$.

Now we bring in the information given
by Lemma \ref{adapt}.
Since we can automatically set some of
the terms in our sequence to zero, it suffices
to prove our result for any subsequence of $D$
we like.
So, we can assume without loss of
generality that
\begin{equation}
|[Pd_k]_2|<10^{-k}.
\end{equation}
Passing to a further subsequence, we can
reduce to the case where either
$[Pd_k]_2$ is always positive or always
negative.  We will consider the positive
case.  The negative case has the same
treatment, except for a sign.
Let $\delta_k=[Pd_k]_2$.  We have
$\delta_k \in (0,10^{-k})$.

Looking at Equation \ref{geodesic},
we see that
$\Xi(Y_D)$ develops into a bounded
subset of the plane.  More precisely,
$\Xi(Y_D)$ develops into the plane as
a translate of the set $C' \times (-C')$,
where $C'$ is the set of all finite sums
of the form
\begin{equation}
\sum_{i=1}^n \epsilon_k \delta_k,
\hskip 30 pt
\epsilon_k \in \{0,1\}.
\end{equation}
But $C'$ is dense in a Cantor set which
consists of all the infinite sums of same
form.  Hence the set of accumulation points
of $\Xi(Y_{\mu})$ contains a Cantor set,
all of whose points belong to $\gamma_P$.
\endproof

\subsection{The End of the Proof}
\label{endproof}

Let $P=2A/(1+A)$ as usual.  Let
$\cal C$ be the set of centers of the
square tiles.
There exists some binary sequence $B$ and some
associated limit realization $\Phi_{\infty}$ such
that
\begin{itemize}
\item The corresponding classifying map
$\Xi_{P,V}$ carries every point of
$\cal C$ into the
interior of a partition piece.  That is, $V$ is a
good offset.  Let $\Xi_{\infty}=\Xi_{P,V}$.
This map is the limit of
maps $\Xi_n-(0,y_n)$. Here
$y_n$ is such that
$y_n+(0,1/2) \in \langle \Gamma_n \rangle$.
\item The curve $\Gamma_{\infty}$
associated to $\Phi_{\infty}$ has
unbounded projections in both coordinate directions.
\item $\Gamma_{\infty}$ is the Hausdorff limit of
the curves $\gamma_n-(0,y_n)$, and this is the same
as the Hausdorff limit of the loops
$\Gamma_n-(0,y_n)$.   The portion of
$(\Gamma_n-\gamma_n)-(0,y_n)$ exits
every compact set in the plane and is not
seen in the limit.
\end{itemize}

Let $c_0=(1/2,1/2)$.  This point
belongs to $\Gamma_{\infty}$.
To finish the proof, we have to check that
$\Gamma_{\infty}$ realizes the
vector dynamics for $\Xi_{\infty}$ associated
to the point $c_0$.  Let $B_R$ denote the
ball of radius $R$ about the origin.

Note that $\Gamma_n-(0,y_n)$ describes
the vector dynamics of $\Xi_n-(0,y_n)$ associated
to $c_0$.
Since $\Xi_{\infty}$ is defined on all of
$\cal C$, we get the following result by continuity.
For any $R>0$ there is some $N$ such that
$n>N$ implies that
$\Xi_n-(0,y_n)$ and $\Xi_{\infty}$ have the same
vector dynamics for relative to all points
of ${\cal C\/} \cap B_R$.
That is, the partition labels assigned to
points of ${\cal C\/} \cap B_R$ by the two
maps are the same.
Hence, $\Gamma_{\infty}$ describes the
vector dynamics of $\Xi_{\infty}$ associated
to $c_0$ on arbitrarily
large balls about the origin.
This is the same as saying that
$\Gamma_{\infty}$ describes
the vector dynamics of $\Xi_{\infty}$ associated
to $c_0$.

This completes the proof of the Unbounded
Orbits Theorem modulo the results quoted
in \S \ref{strong}, the Box Lemma,
and the Copy Theorem.

\newpage
\section{Some Elementary Number Theory}

In this chapter we will justify all the statements
quoted in \S \ref{strong}.  We will also prove
a number of other technical results which
establish useful identities and inequalities
between the rationals discussed in \S \ref{strong}.

\subsection{All About Predecessors}

We fix an even rational parameter
$p/q$.  To avoid trivialities, we assume
that $p>1$.  Given an even rational
parameter $p/q$, we make the
following definitions:
\begin{itemize}
\item $p'/q'$ is the even predecessor of $p/q$.
\item $\widehat p/\widehat q$ is the core
predecessor of $p/q$.
\item $\omega=p+q$.
\item $\tau$ is the tune of $p/q$, as in \S \ref{tune}.
\item $\kappa$ is as in Lemma \ref{strong}.
\end{itemize}
We mean to define $\omega,\tau,\kappa$ for the
other parameters as well. Thus, for instance,
$\omega'=p'+q'$.

\begin{lemma}
\label{omnibus}
The following is true.
\begin{enumerate}
\item $\widehat p/\widehat q \in (0,1)$
is a rational in lowest terms.
\item Either
$\tau-\tau'=\kappa \omega'$ or
$\tau+\tau'=(1+\kappa) \omega'$.
In all cases, $\tau' \leq \tau$.
\item $p'/q'$ is also the even predecessor of
$\widehat p/\widehat q$.
\item $\widehat \kappa=0$.
\item $\omega-2\tau = \widehat \omega-2\widehat \tau$.
\item If $\kappa=0$ then $\tau'=\tau$ when
$\tau<\omega/4$ and $\tau'=\omega'-\tau$ when $\tau>\omega/4$.
\item If $\kappa \geq 1$ then $\kappa \widehat \omega<(3/2)\omega$.
\end{enumerate}
\end{lemma}

\subsubsection{Statement 1}

We first give a formula for $p'/q'$.
 There is some
integer $\theta>0$ so that
$$2p\tau = \theta(p+q) \pm 1.$$
Rearranging this, we get
$$p(2\tau-\theta) - q(\theta) = \pm 1$$
We have $\theta<\min(p,2\tau)$.  Setting 
$$p''=\theta, \hskip 20 pt q''=2\tau-\theta,$$
we see that $|p''q-q''p|=1$.
This implies that $p''/q''$ is in lowest
terms.  Moreover, both $p''$ and $q''$ are
odd.  That means that $p'=p-p''$ and $q'=q-q''$.
Hence
\begin{equation}
\label{zoop}
\omega'=\omega-2\tau.
\end{equation}

We have
\begin{equation}
\label{gap}
\frac{\omega'}{\omega}=1-2\bigg(\frac{\tau}{\omega}\bigg)
\leq 1-2\bigg(\frac{\kappa}{2\kappa+1}\bigg)=\frac{1}{2\kappa+1}.
\end{equation}
If $2\kappa q' \geq q$ then Equation \ref{gap} implies
that $(2\kappa+1)p'+q' \leq p$. 
But then a short calculation shows 
that $|pq'-qp'|=p'q'>1$.
This is a contradiction.
Hence $2\kappa q'<q$.  The same
argument works with $p$ in place of $q$.

Since $\widehat p/\widehat q$ is Farey related
to $p'/q' \in (0,1)$, we see that $\widehat p/\widehat q$
also lies in $(0,1)$ and is in lowest terms.
This proves Statement 1.

\subsubsection{Statement 2}

We first prove that
$\tau \equiv \pm \tau'$ mod $\omega'$.
We will assume that $2p\tau \equiv 1$ mod $\omega$.
The other case has a similar treatment.  We
use the formulas that we have just derived.
We have 
$$2p\tau=\theta \omega +1=
\theta \omega'+2\tau \theta+ 1.$$
Rearranging, we get
$$\theta \omega ' + 1 =
2(p-\theta)\tau = 2p'\tau.$$
Hence $2p'\tau \equiv \pm 2p'\tau'$ mod $\omega'$.
Since $2p'$ is relatively prime to $\omega'$, we
get that $\tau \equiv \pm \tau'$ mod $\omega'$.

Suppose first that
$\tau \equiv \tau'$ mod $\omega'$.
The quantity $t=\tau-\kappa \omega'$ satisfies
$2p't \equiv \pm 1$ mod $\omega'$.
So, we just have to show that
$t \in (0,\omega')$.
We compute
\begin{equation}
t=\tau-\kappa \omega'=\tau - \kappa(\omega-2\tau)=
\tau(2\kappa+1)-\kappa\omega.
\end{equation}
After a bit of algebra, we see that
this expression lies in $(0,\omega')$.

When $\tau \equiv -\tau$ mod $\omega$, we
instead consider the quantity
$t=(\kappa+1)\omega-\tau$.  Again we
know that $2p't \equiv \pm 1$ mod $\omega'$.
The rest of the proof is very similar
to what is done in the other case.

Now we establish the second part of
Statement 1. We have one of the two equalities
$$\tau-\tau'=\kappa \omega', \hskip 30 pt
\tau+\tau'=(\kappa+1)\omega'.$$
In the first case, we obviously
have $\tau' \leq \tau$.  In
the second case, the same
inequality follows from
the fact that $\tau'<\omega'/2$.

\subsubsection{Statement 3}

Let $\widehat p/\widehat q$ be the core
predecessor of $p/q$.
Since
$\widehat p=p-2\kappa p'$ and $\widehat q=q-2\kappa q'$
the two rationals $p'/q'$ and $\widehat p/\widehat q$ are
both Farey related.
Equation \ref{gap}, which is a strict
inequality when $p>1$,  says that
$(2\kappa +1) \omega'< \omega$.
Since $\widehat \omega=\omega-2\kappa \omega'$,
we conclude that
$\omega'<\widehat \omega$.
Hence $p'/q'$ is the Farey predecessor
of $\widehat p/\widehat q$.
This completes the proof of Statement 3.

\subsubsection{Statement 4}

If Statement 4 is false, then
Equation \ref{gap}, applied to
$\widehat p/\widehat q$, tells us that
\begin{equation}
\widehat \omega-3\omega'>0.
\end{equation}
Combining this with the 
fact that 
$$\widehat \omega=\omega-2\kappa \omega'>0,$$
we can say that
$$\omega-(2\kappa+3)\omega'>0.$$
But then we have
$$
\frac{1}{2\kappa+3}=
1-2F_{\kappa+1}<
1-2\bigg(\frac{\tau}{\omega}\bigg)=
\frac{\omega'}{\omega}<
\frac{1}{2\kappa+3}.
$$
This is a contradiction.

\subsubsection{Statement 5}

The proof goes by induction on $\kappa$. 
The result is trivial for $\kappa=0$, so
we assume that $\kappa \geq 1$.
We will change notation somewhat.
Let $p_2/q_2=p/q$ and let
$p_1/q_1=p'/q'$, the even predecessor
of $p_2/q_2$.  Define
\begin{equation}
p_3=p_2-2p_1, \hskip 30 pt
q_3=q_2-2q_1.
\end{equation}
Then
$\kappa_3=\kappa_2-1$.  It suffices to
prove that
$\omega_2-2\tau_2=\omega_3-2\tau_3$.
This is what we will do.

We have
$$w_3=\omega_2-4\tau_2.$$
We define $u_3=3\tau_2-w_2$.  Note that
$\omega_2-2\tau_2=\omega_3-2u_3$.
So, to finish the proof, we just have
to check that $u_3=\tau_3$.  By hypothesis,
we have $u_3>0$.  Since 
$0<\omega_2-2\tau_2=\omega_3-2u_3$, we have
$u_3 \in (0,\omega_3/2)$.  So, we just have
to check that $2p_3u_3 \equiv \pm 1$ mod $\omega_3$.
We will consider the case when $2p_2\tau_2 \equiv 1$ mod $\omega_2$.
The other case works the same way.

In the case at hand, we have
$$k_2=\frac{2p_2\tau_2-1}{\omega_2}.$$
Using this equation, we compute
$$N:=\frac{2p_3u_3-1}{\omega_3}=
\frac{-3-2p_2^2-2p_2q_2 + 6p_2\tau_2}{\omega_2}.$$
Plugging in the equation
$6p_2\tau_2=3+3k_2\omega_2$ and factoring, we
see that $N$ is an integer.  Hence
$2p_3u_3 \equiv 1$ mod $\omega_3$.

\subsubsection{Statement 6}

When $\kappa=0$,
Statement 1
says that either
$\tau=\tau'$ or $\tau+\tau'=\omega'$.
Also, from Equation \ref{zoop} we see that
$2\omega'<\omega$ if and only if
$\tau>\omega/4$.

Suppose that
$\tau'=\tau$ and
$\tau>\omega/4$.
We have $$\tau'=\tau>\omega/4>\omega'/2.$$ which 
contradicts the fact that $\tau'<\omega'/2$.
Suppose that
$\tau'=\omega'-\tau$ and
$\tau<\omega/4$.
We have
$$2\tau'=2\omega'-2\tau>\omega-2\tau=\omega,$$
which again is impossible.

\subsubsection{Statement 7}

We have
$$\widehat \omega=\omega-2\kappa \omega'=
\omega - 2\kappa(\omega-2\tau)=
(1-2\kappa)\omega + 4\kappa \tau.$$
Using this equality, we see that
Statement 7 is equivaent to the truth of
\begin{equation}
\label{refo2}
\frac{\tau}{\omega}<\frac{4\kappa^2-2\kappa+3}{8\kappa^2}.
\end{equation}
By definition
$$
\frac{\tau}{\omega}<\frac{\kappa+1}{2\kappa+3}.
$$
Moreover,
$$\bigg(\frac{4\kappa^2-2\kappa+3}{8\kappa^2}\bigg)
-\bigg(\frac{\kappa+1}{2\kappa+3}\bigg)=\frac{9}{8\kappa^2(2+2\kappa)}>0.
$$
Therefore, Equation \ref{refo2} holds.

\subsection{Existence of the Predecessor Sequence}

We begin by proving the existence of
an auxiliary sequence.  The notation
$p'/q' \leftarrow p/q$ means that
$p'/q'$ is the even predecessor of $p/q$.

\begin{lemma}
\label{weak}
Let $A \in (0,1)$ be irrational.
There exists a sequence
$\{p_n/q_n\}$ converging to
$A$ such that
$p_n/q_n \leftarrow p_{n+1}/q_{n+1}$
for all $n$.
\end{lemma}

\startproof
Let $\H^2$ denote the upper half-plane model
of the hyperbolic plane.
We have the usual Farey
triangulation of $\H^2$ by
ideal triangles.  The geodesics
bounding these triangles join
rationals $p_1/q_1$ and $p_2/q_2$
such that $|p_1q_2-p_2q_1|=1$.
We call these geodesics the
{\it Farey Geodesics\/}.
Two of the Farey geodesics in the tiling
join the points $0/1$ and $1/1$ to $1/0$.
This last point is interpreted as the
point at infinity in the upper half
plane model of $\H^2$.  Let
$Y=\{(x,y)|\ 0<x<1\}$ be the
open portion of $\H^2$ between these
two vertical geodesics.

Each Farey geodesic in $Y$ which joins
two even rationals $p/q$ and $p'/q'$ 
can be oriented
so that it points from $p/q$ to $p'/q'$ if
and only if $p'/q' \leftarrow p/q$.
We leave the rest of the Farey geodesics
unoriented.

We claim that there is a backwards 
oriented path from
$0/1$ to $A$.  To see this, choose any
sequence of even rationals converging to $A$
and consider their sequence of even predecessors.
This gives us a sequence of finite directed
paths joining $0/1$ to rationals which
converge to $A$.  
Given the local finiteness
of the Farey triangulation -- meaning that
any compact subset of $\H^2$ intersects
only finitely many triangles -- we can
take a limit of these paths, at least on
a subsequence, and get a directed path
converging to $A$.  

Reading off the
vertices of this path gives us a
sequence $\{p_n/q_n\}$ with
$p_n/q_n \leftarrow p_{n+1}/q_{n+1}$
for all $n$.
\endproof

\noindent
{\bf Proof of Lemma \ref{neat}:\/}
We get the predecessor sequence from the
sequence in Lemma \ref{weak} in
the following way.  Before each
term $p_n/q_n$ with $\kappa_n \geq 1$
we insert the rational $\widehat p_n/\widehat q_n$.
This is the core predecessor of $p_n/q_n$.
We know from the lemma above that
$\widehat \kappa_n=0$ and that
$p_{n-1}/q_{n-1}$ is the predecessor
of $\widehat p_n/\widehat q_n$. In short
$$p_{n-1}/q_{n-1} \prec \widehat p_n/\widehat q_n \prec p_n/q_n.$$
Once we make all these insertions, the
resulting sequence is a predecessor sequence which
converges to $A$.  Again, we will not stop
to prove uniqueness because we don't care
about it.
\endproof

\subsection{Existence of the Approximating Sequence}

We call the sequence from Lemma \ref{weak} the
{\it even predecessor sequence\/}.  As we
have just remarked, the even predecessor
sequence is a subsequence of the predecessor sequence.

\begin{lemma}
\label{neat1}
Let $\{p_n/q_n\}$ be the
even predecessor sequence
which converges to $A$.
There are infinitely many
values of $n$ such that
$\tau_{n+1}>\omega_{n+1}/4$.
\end{lemma}

\startproof
If this lemma is false, then we might
chop off the beginning and assume that
this never happens. 
We have the formula 
$\omega_n=\omega_{n+1}-2\tau_{n+1}$.
So, if this lemma is false then
we have $2\omega_n<\omega_{n+1}$ only
finitely often.   Chopping off
the beginning of the sequence, we
can assume that this never occurs.

We never have $2\omega_n=\omega_{n+1}$
because $\omega_{n+1}$ is odd.  So, we
always have $\omega_n>\omega_{n+1}$.
In this case, we
have $$r_{n-1}=2r_n-r_{n+1}$$ for
$r \in \{p,q\}$.  Applying this iteratively,
we get that $r_{k}=(k+1)r_1-kr_0$ for $k=1,2,3,...$.
But then $\lim p_n/q_n=(p_1-p_0)/(q_1-q_0)$, which
is rational. This is a contradiction.

Now let $p_n/q_n$ be a term with
$2\omega_n<\omega_{n+1}$.
We will consider the case when $p_n/q_n<p_{n+1}/q_{n+1}$.
The other case has a similar treatment.
We introduce the new rational
$$p_n^{*}/q_n^{*}=\frac{p_{n+1}-p_n}{q_{n+1}-q_n}.$$
Note that the three rationals
$p_n/q_n$ and $p_{n+1}/q_{n+1}$ and
$p_n^{*}/q_n^{*}$ form the vertices of a
Farey triangle.  Moreover $p_{n+1}/q_{n+1}<p_n^*/q_n^*$
because $p_{n+1}/q_{n+1}$ is the Farey sum of
$p_n/q_n$ and $p_n^*/q_n^*$.

Given the geometry of the Farey graph, we have
$$A \in [p_n/q_n,p_n^*/q_n^*].$$
Moreover
$$|p_n/q_n-p_n^*/q_n^*|=\frac{1}{q_nq^*_n} \leq
\frac{1}{q_n^2}<\frac{4}{\omega_n^2}.$$
The second-to-last inequality comes from the
fact that $q_n^*>q_n$.
\endproof

\noindent
{\bf Proof of Lemma \ref{strong0}:\/}
Let $\{p_k/q_k\}$ be the even predecessor
sequence.  If it happens infinitely often that
$\kappa_k \geq 1$ then there are infinitely
many core terms in the predecessor sequence.
Otherwise
the predecessor sequence and the even predecessor
sequence have the same tail end.  But then
the previous result shows that there are
infinitely many values of $k$ for which
$2\omega_k<\omega_{k+1}$.  This gives
infinitely many strong terms.
Suppose that $p/q, p^*/q^*$ are two
consecutive terms in the predecessor sequence,
with $p/q$ non-weak.
\newline
\newline
{\bf Case 1:\/}
Suppose that $p/q$ and $p^*/q^*$ are
both terms in the even sequence.
Then $p/q \leftarrow p^*/q^*$ and
$\tau_*>\omega_*/4$.  But then
$2\omega<\omega_*$, as in 
the proof of Lemma \ref{neat1}.
From here, we get the Diophantine estimate
from Lemma \ref{neat1}.
\newline
\newline
{\bf Case 2:\/}
Suppose that neither $p/q$ not $p^*/q^*$
are terms in the even sequence.  This
case cannot happen because when
$p^*/q^*$ is the core predecessor of
some rational and $p/q$ is the core
predecessor of $p^*/q^*$.  This
contradicts Statement 4 of Lemma \ref{omnibus}.
\newline
\newline
{\bf Case 3:\/}
Suppose that $p/q$ is a term of the
even sequence but $p^*/q^*$ is not.
Then the term in the predecessor sequence
after $p^*/q^*$ is $p^{**}/q^{**}$,
which belongs to the even sequence.
By Statement 3 of
Lemma \ref{omnibus}, we have
$p/q \leftarrow p^{**}/q^{**}$.
Moreover $2\omega<\omega^*<\omega^{**}$.
From here, we get the Diophantine estimate
from Lemma \ref{neat1}.
\newline
\newline
{\bf Case 4:\/}
Suppose that $p/q$ is not a term in
the even sequence but $p^*/q^*$ is.
Then $p/q$ is the core predecessor of $p^*/q^*$.
Statement 4 of Lemma \ref{omnibus} says
that $\kappa=0$.  Hence the term
$p'/q'$ preceding $p/q$ in the
predecessor sequence is the even predecessor of $p/q$.
Since $\kappa=0$, we have
$\tau<\omega/3$.  Using the formula
$\omega'=\omega-2\tau$, we see that
$3\omega'>\omega$.  At the same time
$\omega'<\omega<2\omega^*$.  Lemma
\ref{neat1} gives us
$$\bigg|A-\frac{p'}{q'}\bigg|<
\frac{4}{(\omega')^2}<\frac{36}{\omega^2}.$$
At the same time
$$\bigg|\frac{p}{q}-\frac{p'}{q'}\bigg| =
\frac{1}{qq'}<\frac{4}{\omega \omega'}<\frac{12}{\omega^2}.$$
The result in this case then follows from
the triangle inequality.
\endproof

\subsection{Another Identity}
\label{hard}

Let $p/q$ be an even rational parameter with
$\kappa \geq 1$. Let
$\widehat p/\widehat q$ be the core predecessor.
Define
\begin{equation}
w=\widehat \tau, \hskip 30 pt
h=\omega-2\tau.
\end{equation}
Here is an identity between these
quantities.

\begin{lemma}
\begin{equation}
\label{MAIN}
(2\kappa+1) h+2w=\omega.
\end{equation}
\end{lemma}

\startproof
Define
\begin{equation}
F_{\kappa}=\frac{\kappa}{2\kappa+1}.
\end{equation}
Let $F$ be the union of all such
rationals.

We rescale the first block so that it
coincides with the unit square.
This rescaling amounts to dividing by
$\omega$. For each quantity $\lambda$
defined above,
we let $\lambda^*=\lambda/\omega$.
Within the unit square, our identity
is
\begin{equation}
\label{id2}
(2\kappa+1) h^*+2w^*=1.
\end{equation}
The rescaled horizontal and
vertical lines depend continuously
on the parameter $\tau^*$.
Moreover, when $\tau^* \in (F_{\kappa},F_{\kappa+1})$
is rational,
all the horizontal and vertical lines of
capacity up to $4\kappa+2$ are distinct,
because the denominator of the corresponding
rational is at least $2\kappa+1$.
Therefore, $\kappa$ is a locally
constant function of $\tau^*$, and
changes only when $\tau^*$
passes through a value of $F$.

Consider what happens when $\tau^*=F_{\kappa}+\epsilon$,
and $\epsilon \in (0,F_{\kappa+1}-F_{\kappa})$.
We compute
\begin{equation}
\label{derive2}
(2\kappa +1)\tau^* - \kappa = (2\kappa +1) \epsilon<1/2
\end{equation}
Hence $w=(2\kappa+1)\epsilon$, which means that
$2w^* = (4 \kappa +2) \epsilon$.
We also have
$h^*=1-2\tau^*$.
From these equations one can easily compute that
Equation \ref{id2} holds for
all $\epsilon \in (0,F_{\kappa+1}-F_{\kappa})$.
\endproof

\begin{lemma}
\label{right}
$w=\widehat \tau$.  In other words,
the line $x=\widehat \tau$ has
capacity $4\kappa + 2$ with
respect to $p/q$.
\end{lemma}

\startproof
In view of Equation \ref{MAIN}, we just
have to prove that
\begin{equation}
\label{goal5}
(2\kappa+1)(\omega-2\tau) + 2\widehat \tau = \omega
\end{equation}
By Statement 5 of Lemma \ref{omnibus} we have
$\omega-2\tau = \widehat \omega-2\widehat \tau$.
But then Equation \ref{goal5} is equivalent to
\begin{equation}
\label{goal6}
(2\kappa+1)\widehat \omega - 4\kappa \tau = \omega.
\end{equation}
This equation holds, because
$$\omega-\widehat \omega=
2\kappa \omega'=2\kappa(\omega-2\tau)=
2\kappa(\widehat \omega-2\widehat \tau)=
2\kappa \widehat \omega-4\kappa \widehat \tau.$$
This completes the proof.
\endproof

\newpage
\section{The Box Lemma and The Copy Theorem}

\subsection{The Weak and Strong Copying Lemmas}

Let $p/q$ be an even rational parameter
and let $p'/q'$ be the even predecessor of
$p/q$. 
Let $\Pi$ and $\Pi'$ denote the
plaid tilings with respect to these two
parameters. In this chapter we will consider
the case when $\kappa=0$. 
There are two subcases:
\begin{itemize}
\item {\bf weak:\/} $2\omega'>\omega$.
This corresponds to $\tau<\omega/4$.
\item {\bf strong:\/} $2\omega'<\omega$.
This corresponds to $\tau>\omega/4$
(and $\tau<\omega/3$.)
\end{itemize}

Let $R_{p'/q'}$ be the rectangle associated
to the parameter $p'/q'$ in \S \ref{mb}.  Again,
this rectangle is bounded by the lines
\begin{itemize}
\item $y=0$.
\item $y=\omega'$.
\item $x=0$.
\item $x=\min(\tau',\omega-2\tau')$.
\end{itemize}
The right side of $R'_{p'/q'}$ is whichever
line of capacity at most $4$ is closest to
the $y$-axis.

In the weak case, we
let $\Sigma'$ denote the subset of
of $R_{p'/q'}$ bounded above by the line
\begin{equation}
y=\omega'-\min(\tau',\omega'-2\tau').
\end{equation}
That is, the top of $\Sigma'$ is whichever
horizontal line of capacity at most
$4$ lies closest to the
top of $R_{p'/q}$.

In the strong case, we
let $\Sigma'=R_{p'/q'}$.

In both cases we define
\begin{equation}
\Sigma=\Sigma'.
\end{equation}
Even though $\Sigma=\Sigma'$ it is useful
to have separate notation, so that in
general an object $X'$ corresponds to
the parameter $p'/q'$ and an object $X$
corresponds to the parameter $p/q$.

In the next chapter we will prove the following result.
\begin{lemma}[Weak and Strong Copying]
$\Sigma' \cap \Pi'=\Sigma \cap \Pi$.
\end{lemma}

\subsection{The Core Copy Lemma}

Let $p/q$ be an even rational parameter
with $\kappa \geq 1$. Let
$\widehat p/\widehat q$ be the core
predecessor of $p/q$.
Let $\Pi$ and $\widehat \Pi$ denote the
plaid tilings with respect to these two
parameters.

Let $\Upsilon$ be vertical translation by
$(\omega+\widehat \omega)/2$.  It follows
from Statement 5 of Lemma \ref{omnibus} that
$\Upsilon$ maps the horizontal
lines of capacity $\pm 1$ w.r.t
$\widehat p/\widehat q$ to the lines of
capacity $\pm 1$ w.r.t. $p/q$.

Let $R_{\widehat p/\widehat q}$ is the rectangle associated
to the parameter $\widehat p/\widehat q$ in \S \ref{mb}.
Define
\begin{equation}
\widehat \Sigma=R_{\widehat p/\widehat q},
\hskip 30 pt \Sigma=\Upsilon(\widehat \Sigma).
\end{equation}

In \S 9 we will prove the following result.

\begin{lemma}[Core Copying]
$\Sigma \cap \Pi=
\Upsilon(\widehat \Sigma \cap \widehat \Pi)$.
\end{lemma}

\subsection{Proof of the Box Lemma}

Our proof goes by induction on the denominator
of the parameter.  We will suppose that $p/q$ is
a rational with the smallest denominator for
which we don't know the truth of the result.

By construction, $\Gamma$ can only intersect
$\partial R$ in the right edge, which we call $\rho$.
If we knew that $\Gamma$ intersects $\rho$ twice,
then, because $\Gamma$ is a closed loop,
we could conclude that
the portion of $\Gamma$ contained in $\Sigma$ is
an arc whose endpoints are on $\rho$.
If $\rho$ has capacity $2$, then
$\Gamma$ can intersect $\rho$ at most twice,
and the large $x$-diameter of $\Gamma$
implie that $\Gamma$ does intersect $\rho$.
By symmetry, $\Gamma$ intersects
$\rho$ twice. 

Now suppose that $\rho$ has capacity $4$.
In this situation, we have $\kappa \geq 1$,
so we can apply the Core Copy Lemma to
$p/q$ and $\widehat p/\widehat q$.

\begin{lemma}
\label{nocross}
$\Sigma \subset R$.
\end{lemma}

\startproof
The left edges of $\Sigma$ and $R$ both
lie in the $x$-axis.  So, we just have
to show that the width of $\Sigma$ is
at most the width of $R$.
The width of $R$ is $\omega-2\tau$.
Statement 2 of Lemma \ref{omnibus} says
that $\widehat \kappa=0$.  Hence,
the width of $\widehat \Sigma$ is
$\widehat \tau$, a quantity which
not greater than $\widehat \omega-2\widehat \tau$.
(The latter quantity is the distance between
the nearest line of capacity $4$ to the
$x$-axis.)  Now, using
Statement 5 of Lemma \ref{omnibus}, we observe that
$\widehat \tau \leq \widehat \omega - 2\widehat \tau=\omega-2\tau.$
\endproof

\begin{lemma}
$\widehat \Gamma$ does not intersect
the right edge of $\widehat \Sigma$.
\end{lemma}

\startproof
Since $\widehat \kappa=0$, the right
edge of $\widehat \Sigma$ has capacity $2$.
$\widehat \Gamma$ only intersects the line
containing this edge twice, and these 
intersection points must outside of
$\widehat \Sigma$, for otherwise
$\widehat \Gamma$ could not make a
closed loop. Compare Figure 2.1.
\endproof

It now follows from the Core Copy Lemma that
$\Gamma$ does not intersect
the right edge of $\Sigma$.
Figure 6.1
shows three sub-rectangles
$R_1,R_2,\Sigma \subset R$.  The horizontal dividers in 
the picture are the horizontal lines of capacity
$2$.  We have already shown that
$\Gamma$ does not intersect the right
edge of $\Sigma$.  But then $\Gamma$ cannot
cross $\rho$ between the two horizontal
lines of capacity $2$.  To finish our proof,
we just have to check that $\Gamma$ intersects
the right edge of $R_1$ once.  By symmetry,
$\Gamma$ intersects the right edge of $R_2$ once,
and the right edge of $\Sigma$ blocks $\Gamma$
so that it cannot intersect $\rho$ anywhere else.

\begin{center}
\resizebox{!}{2.5in}{\includegraphics{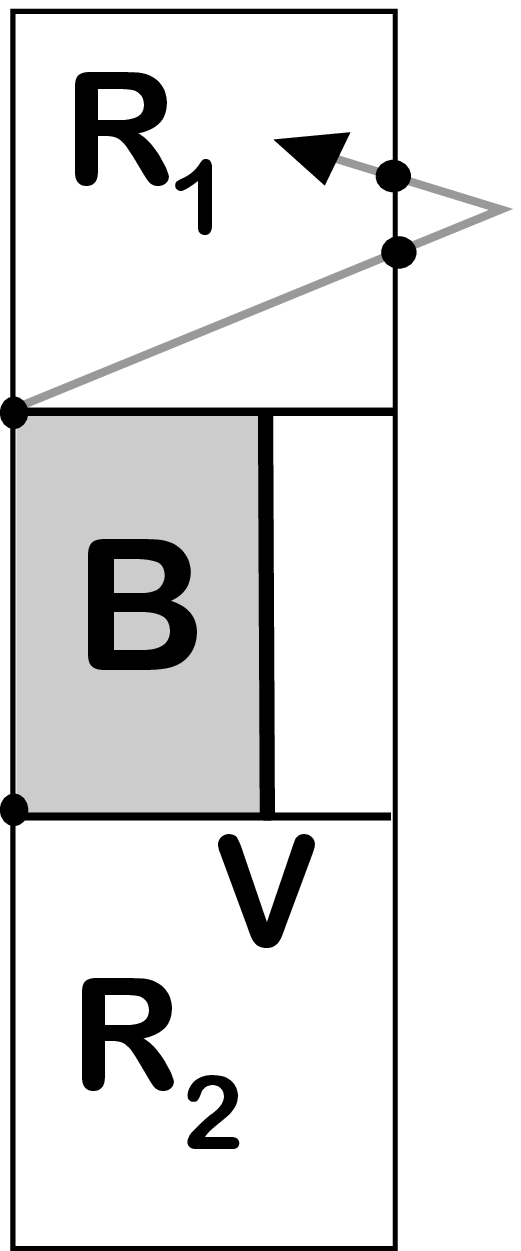}}
\newline
{\bf Figure 6.1:\/} The rectangles $R_1,S,R_2$.
\end{center} 

We know that $\Gamma$ cannot intersect
$\rho_1$ more than twice because then,
by symmetry, $\Gamma$ would intersect
$\rho$ at least $6$ times.
If $\Gamma$ intersects $\rho_1$ exactly
twice, then $\Gamma$ is trapped in
$R_1$ and cannot get around to close up
with the portion of $\Gamma$ outside $R_1$.
Since $\Gamma$ cannot get trapped in this way,
we see that $\Gamma$ intersects $\rho_1$
at most once.  On the other hand,
if $\Gamma$ does not intersect
$\rho_1$ at all, then $\Gamma$ is
trapped in the same way. In short,
$\Gamma$ intersects $\rho_1$ exactly
once, and the same goes for $\rho_2$.

This completes the proof of 
the Box Lemma.

\subsection{Setup for the Copying Theorem}

Let $R_{p/q}$ be the rectangle associated
to the parameter $p/q$.

\begin{lemma}
\label{width}
Let $p'/q'$ be the even predecessor of $p/q$.
Suppose that $\kappa=0$.  Then
$R' \subset R$.
\end{lemma}

\startproof
Both boxes contain $(0,0)$ as the bottom left
vertex.  So, we just have to show that the width
and height of $R'$ are at most the width and
height of $R$.
The width of $R$ is $\tau$.
The width of $R'$ is either
$\tau'$ or $\omega'-2\tau'$,
whichever is smaller.
In either case, we have
${\rm width\/}(R') \leq \tau' \leq \tau={\rm width\/}(R).$
The height of $R$ is $\omega$ and
The height of $R'$ is $\omega-2\tau$.
\endproof

Suppose that $p/q$ and $p^*/q^*$ are two
consecutive terms in the approximating sequence.
We write
\begin{equation}
p/q = p_0/q_0 \prec p_1/q_1 \leftarrow \cdots 
\leftarrow p_n/q_n.
\end{equation}
If $p/q$ is strong then $p_0/q_0$ is the
even predecessor of $p_1/q_1$, and
$p_k/q_k$ is weak for $k=1,...,n$.
If $p/q$ is core, then $p_0/q_0$ is the
core predecessor of $p_1/q_1$ and
$p_1/q_1$ is non-core, and
$p_k/q_k$ is weak for $k=2,...,n$.
Here $k \geq 2$.  

We introduce some notation to help with the
proof.  Let $(\Sigma_k',\Sigma_{k+1})$ be the
pair of rectangles associated to the
pair of parameters $(p_k/q_k,p_{k+1}/q_{k+1})$.
Note that
$p_k/q_k$ gets two such rectangles attached
to it, namely $\Sigma_k$ and $\Sigma'_k$.
These rectangles play different roles in
the proof. Let $\Pi_k$ denote the plaid
tiling associated to $p_k/q_k$.
In general, we set $\omega_k=\omega(p_k/q_k)$, etc.

\subsection{The Strong Case}

Suppose first that $p_0/q_0$ is strong.  
Recalling that $\Sigma_0'=\Sigma_1=B_0$, we conclude that
\begin{equation}
B_0 \cap \Pi_0=\Sigma_0' \cap \Pi_0=^*\Sigma_1 \cap \Pi_1=
B_0 \cap \Pi_1.
\end{equation}
The starred equality comes from the
Strong Copy Lemma.

Recall that $TH_k$ and $BH_k$ are the
top and bottom horizontal lines of
capacity $2$ with respect to $p_k/q_k$.

\begin{lemma}
\label{linematch}
For each $k=1,...,n$ the rectangle
$B_0$ is contained in the lower half of $B_k$
and one of $BH_0$ or $TH_0$ coincides
with $BH_k$.
\end{lemma}

\startproof
Note that every box in sight contains
$(0,0)$ as the lower left vertices.
So, we can decide which box contains
which other box just by looking at the
widths and heights.  

Eince $\tau_1>1/4$ we have $2\omega_0<\omega_1$,
the height of $B_0$ is less than half that
of $B_1$.  Hence $B_0$ lies in the
lower half of $B_1$.   By Lemma
\ref{width}, we have
$B_1 \subset ... \subset B_k$.
Hence $B_0$ lies in the lower
half of $B_k$.

From Statement 1 of Lemma \ref{omnibus},
we have either $BH_0=BH_1$ or $TH_0=BH_1$.
From Statement 6 of Lemma \ref{omnibus}, we
see that $BH_1=...=BH_k$.
Hence either $BH_0=BH_k$ or
$TH_0=BH_k$.
\endproof

Note that $\Sigma_k'$ contains
the lower half of $B_k$.  Hence
$B_0 \subset \Sigma_k'$ for all $k$.
Now we will show inductively that
$B_0 \cap \Pi_k$ implies
$B_0 \cap \Pi_{k+1}$.
We already have proved this for $k=0$.
For each $k=1,...,n-1$, we have
$\tau_k<1/4$.  We have

\begin{equation}
B_0 \cap \Pi_0=B_0 \cap \Pi_k \subset
\Sigma'_k \cap \Pi_k=^*\Sigma_{k+1} \cap \Pi_{k+1}.
\end{equation}
The starred equality comes from the Weak
Copy Lemma.
From this equation, we get
$B_0 \cap \Pi_0=B_0 \cap \Pi_{k+1}$.
By induction,
$B_0 \cap \Pi_0=B_0 \cap \Pi_n$.  
That is, $\Pi_n$ copies $\Pi_0$ inside $B_0$, and
$B_0$ lies in the lower half of $B_n$.
This is what we wanted to prove.

\subsection{The Core Case}

Suppose that $p_0/q_0$ is core.
Let $\Upsilon$ be the vertical translation by
$(\omega_1-\omega_2)/2$.  
By construction $\Upsilon(B_0)$ is symmetric
with respect to the horizontal midline of $B_1$.
Statement 5 of Lemma \ref{omnibus} says that
the distance between the two horizontal lines
of capacity $2$ is the same w.r.t $p_0/q_0$ and
w.r.t. $p_1/q_1$. Hence, by symmetry,
$\Upsilon(TH_0)=TH_1$ and
$\Upsilon(BH_0)=BH_1$.
We have
\begin{equation}
\Upsilon(B_0 \cap \Pi_0)=
\Upsilon(\Sigma_0' \cap \Pi_0)=\Sigma_1 \cap \Pi_1=
\Upsilon(B_0) \cap \Pi_1.
\end{equation}

\begin{lemma}
$\Upsilon(B_0) \subset \Sigma_1'$.
\end{lemma}

\startproof
The left edge of $\Upsilon(B_0)$ lies in the
$y$-axis, just like the left edge of $\Sigma_1'$.

The width of $\Upsilon(B_0)$ is $\tau_0$ and the
width of $\Sigma_1'$ is $\omega_1-2\tau_1$.
By Statement 2 of Lemma \ref{omnibus}, the
predecessor sequence cannot have $2$ core terms in
a row.  Hence $\kappa_0=0$.
This means that
$3\tau_0<\omega_0$.
Hence
$$\tau_0<\omega_0-2\tau_0=\omega_1-2\tau_1.$$
The equality is Statement 5 of Lemma \ref{omnibus}.
This takes care of the widths.

The $y$-coordinate of the top edge of
$\Upsilon(B_0)$ is 
$(\omega_0+\omega_1)/2$.
The height of $\Sigma_1'$ is either
$\omega_1$ or $2\tau_1$ depending
on whether $p_1/q_1$ is strong or weak.
Since $2\tau_1<\omega_1$, we just need
to deal with the weak case.  That is,
we have to show that
$\omega_0+\omega_1<4\tau_1.$
But 
$4\omega_1/3<4\tau_1$
because $p_0/q_0$ is core.
So, it suffices to show that
$$\omega_0<\omega_1/3.$$
But this follows from
Equation \ref{gap}, because
$\kappa_1 \geq 1$.
\endproof

Now we can take the next step in the
argument.  
\begin{equation}
\Upsilon(B_0 \cap \Pi_0)=
\Upsilon(B_0) \cap \Pi_1 
\subset \Sigma_1' \cap \Pi_1 =^*
\Sigma_2 \cap \Pi_2.
\end{equation}
The starred equality is either the Weak Copy
Lemma or the Strong Copy Lemma, whichever applies.
Hence $\Upsilon(B_0 \cap \Pi_0)=\Upsilon(B_0) \cap \Pi_2$.

\begin{lemma}
For each $k=2,...,n$ the rectangle
$\Upsilon(B_0)$ is contained in the lower half of $B_k$
and one of $BH_0$ or $TH_0$ coincides with
$BH_k$.
\end{lemma}

\startproof
By Lemma \ref{width}, the widths
of $R_1,...,R_k$ are non-decreasing.
So, the width of $\Upsilon(B_0)$ is
at most the width of $R_k$.

The $y$-coordinate of the top edge of
$\Upsilon(B_0)$ is
$(\omega_0+\omega_1)/2$.  The height
of $B_2$ is $\omega_2$.
By Statement 8 of Lemma \ref{omnibus}, we see that
the $y$ coordinate of the top edge
of $\Upsilon(B_0)$ is less than half
the height of $B_2$.  This takes
care of the case $k=2$.

By Lemma \ref{width}, we have
$B_2 \subset ... \subset B_k$,
and all these boxes have $(0,0)$
as their bottom left vertex.  Hence,
$\Upsilon(B_0)$ lies in the
bottom half of $B_k$.

We have already seen that
$\Upsilon(BH_0)=BH_1$ and
$\Upsilon(TH_0)=TH_1$.
By Lemma \ref{linematch}, one
of these two lines coincides with $BH_2$.
But then Statement 6 of Lemma \ref{omnibus} says that
$BH_2=...=BH_k$.
\endproof

The rest of the proof is the same inductive argument
as in the Strong Case.

\newpage
\section{The Weak and Strong Copy Lemmas}

\subsection{Geometric Alignment}
\label{align}

In this chapter we prove the Weak and Strong
Copy Lemmas.  The two results have essentially
the same proof, except for a few minor details.
The idea is to verify the conditions of the
Matching Lemma from \S \ref{MATCH}.
We use the notation and terminology
from \S \ref{MATCH}. Here we
have $\Upsilon={\rm Identity\/}$ and
hence $\Sigma=\Sigma'$.

The two plaid tilings agree along the bottom
edge of $\Sigma$ and $\Sigma'$, because this
common edge lies in the boundary of the
first block w.r.t. both parameters.
Hence $(\Sigma,\Pi)$ and $(\Sigma',\Pi')$ are
weakly horizontally aligned.

\begin{lemma}
\label{matrix}
$(\Sigma,\Pi)$ and 
$(\Sigma',\Pi')$ are geometrically aligned.
\end{lemma}

\startproof
Let $z$ and $z'$ be corresponding points
in $\Sigma=\Sigma'$.  These points
lie on slanting lines of the same type
which have the same $y$-intercept.
The difference in slopes of the two lines is
$$|P-P'|=|Q-Q'|=
\frac{2}{\omega\omega'}.$$
Hence
\begin{equation}
\|z-z'\| \leq \frac{2\tau'}{\omega \omega'}<\frac{1}{\omega}<\frac{1}{\omega'}
\end{equation}
But $z'$ is at least
$1/\omega'$ from the interval contaiing it.
Hence $z$ and $z'$ lie in the same vertical unit interval.
\endproof

\noindent
{\bf Remark:\/}
In replacing $1/\omega$ by $1/\omega'$, we threw away
some of the strength of our estimate.  We did this
because in \S 9 we will have a much tighter estimate,
and we want the two arguments to look similar.

\subsection{Alignment of the Capacity Sequences}
\label{capacityX}

We introduce new variables
\begin{equation}
M_i=M_i\omega, \hskip 30 pt
C_j=c_j\omega, \hskip 30 pt
M_i'=m_i'\omega', \hskip 30 pt
C_j'=c_i'\omega'.
\end{equation}

\begin{itemize}
\item $C_i$ is a nonzero even integer in $(-\omega,\omega)$.
\item $C_i'$ is a nonzero even integer in $(-\omega',\omega')$.
\item $M_j'$ is an odd integer in $[-\omega',\omega']$.
\item $M_j$ is an odd integer in $[-\omega,\omega]$.
\end{itemize}

Let 
\begin{equation}
\lambda= \omega'/\omega.
\end{equation}
We have
\begin{equation}
 C_i'=[2 \omega'  P' i]_{2 \omega'},
\hskip 20 pt
C_i=[2\omega P i]_{2 \omega},
\hskip 20 pt 
\lambda C_i=[2 \omega' P i]_{2  \omega'}.
\end{equation}

We introduce the expression
\begin{equation}
\Psi(i)=|2 \omega Pi - 2 \omega'   P' i|.
\end{equation}
As long as
\begin{equation}
\label{key0}
\Psi(i)<\ell + 2\lambda, \hskip 30 pt
\min(| C_i'|, \omega'-|C_i'|) \geq \ell
\end{equation}
the signs of $ C_i'$ and $C_i$ are the same.

Using the fact that
\begin{equation}
|P-P'|=\frac{2}{\omega\omega'}
\end{equation}
we see that
\begin{equation}
\label{bound0}
\Psi(i)=\frac{4 i}{\omega} \leq \frac{4W'}{\omega}.
\end{equation}

In all cases, we have $W' \leq \tau' \leq \omega'/2$.
Hence $$\Psi(i)<2\lambda<1+2\lambda.$$
So, we can take $\ell=1$ in Equation \ref{key0}.
Hence, the two capacity sequences are
arithmetically aligned.

\subsection{Alignment of the Mass Sequences}
\label{masschange1}

We have
\begin{equation}
 M_i'=[\omega'  P' i+\omega']_{2 \omega'},
\hskip 20 pt
M_i=[2\omega P i+\omega]_{2 \omega},
\hskip 20 pt 
\lambda M_i=[ \omega' P i+\omega']_{2  \omega'}.
\end{equation}

We introduce the expression
\begin{equation}
\Xi(i)=|\omega Pi - \omega'   P' i|.
\end{equation}
This function differs from $\Psi$ just in a factor of $2$.
As long as $m_i'$ is a signed term and
\begin{equation}
\label{key1}
\Xi(i)<\ell + \lambda, \hskip 30 pt
\min(| M_i'|, \omega'-|M_i'|) \geq \ell
\end{equation}
the signs of $ M_i'$ and $M_i$ are the same.

We have
\begin{equation}
\label{xi}
\Xi(i)=\frac{2i}{\omega}<\frac{2H'+4W'}{\omega}.
\end{equation}
There are several cases to consider.
\newline
\newline
{\bf Case 1:\/}
Suppose that we are in the weak case and
$\tau' \leq 2\omega'-\tau'$. 
We have 
$$H'=\omega'-\tau', \hskip 30 pt
W'=\tau'.$$ 
Using the facts that
$$\tau'=\tau, \hskip 30 pt
\omega'=\omega-2\tau,$$
and plugging the values of $W'$ and $H'$
into Equation \ref{xi}, we get
$$\Xi(i)<\frac{2\omega'+2\tau'}{\omega}=
\frac{2\omega'+2\tau}{\omega}=
\frac{2\omega-2\tau}{\omega}=
\frac{\omega+\omega'}{\omega}=1+\lambda.$$
So, Equation \ref{key1} holds
and there are no sign changes.
\newline
\newline
{\bf Case 2:\/}
Suppose that we are in the weak case and
$\tau'>2\omega'-\tau'$.   This
means that $\tau<\omega/4$ and
$\tau'=\tau$ and $\tau'>\omega'/3$.
We have 
$$H'=2\tau', \hskip 30 pt
W'=\omega-2\tau'.$$
Plugging this into Equation \ref{xi},
we get
$$\Xi(i)<\frac{4\omega'-4\tau'}{\omega}<\frac{8\tau'}{\omega}=
\frac{8\tau}{\omega}<2.$$
So, we only have to worry about the case
when $\widehat C_i=\pm 1$.  This happens
for $i \in \{-\tau',\tau',\omega'-\tau'\}$.
The highest index is $\omega'+\tau'-1$.
we have

$$
\Psi(\pm \tau') \leq
\Psi(\omega'-\tau')=\frac{2\omega'-2\tau'}{\omega}=
\frac{2\omega'-2\tau}{\omega}=\frac{2\omega-4\tau}{\omega}=$$
\begin{equation}
\frac{\omega+(\omega-4\tau)}{\omega}<\frac{\omega+\omega'}{\omega}=1+\lambda.
\end{equation}
So Equation \ref{key1} holds in all cases.
\newline
\newline
{\bf Case 3:\/}
Suppose we are in the strong case and
$\tau' \leq 2\omega'-\tau'$.
We have 
$$H'=\omega', \hskip 30 pt
W'=\tau'.$$
This gives us the estimate
\begin{equation}
\label{bound2}
\Xi(i)<\frac{2\omega'+4\tau'}{\omega}.
\end{equation}
Since $3\tau' \leq \omega'$ and $\tau+\tau'=\omega'$, 
we have
$$2\tau'\leq \omega'-\tau'=\tau.$$
Hence
$$\Xi(i)<\frac{2\omega'+2\tau}{\omega}=
\frac{(\omega'+2\tau)+\omega'}{\omega}=1+\lambda.$$
So, Equation \ref{key1} holds in all cases,
and there are no sign changes.
\newline
\newline
{\bf Case 4:\/}
Suppose we are in the strong case and
$\omega'-2\tau'<\tau'$.  In this case, we have
$$H'=\omega', \hskip 30 pt
W'=\omega'-2\tau'.$$
This gives us the estimate
\begin{equation}
\label{bound3}
\Xi(i)<\frac{6\omega'-8\tau'}{\omega}.
\end{equation}
Using the fact that $\tau+\tau'=\omega'$
and $\tau<\omega/3$, we get
$$\Xi(i)<\frac{6\tau-2\tau'}{\omega}<2.$$
Once again, we just have to worry about
$C_i=\pm 1$.  The relevant
indices are $i \in \{-\tau',\tau',\omega'-\tau',\omega'+\tau'\}$.
We have
$$\Psi(\pm \tau') \leq \Psi(\omega'-\tau')<
\Psi(\omega'+\tau') =
\frac{2\omega'+2\tau'}{\omega}<
\frac{\omega+\omega'}{\omega}=1+\lambda.$$
Here we used the fact that $2\omega'<\omega$
and $2\tau'<\omega'$.
So, Equation \ref{key1} holds in all
cases, and the mass sequences are aligned.

In short $(\Sigma,\Pi)$ and $(\Sigma',\Pi')$
are arithmetically aligned.  We have
verified all the conditions of the
Matching Lemma, and so 
$\Sigma \cap \Pi=\Sigma'\cap \Pi'$. This
proves the Weak and Strong Copy Lemmas.

\newpage

\section{The Core Copy Lemma}

\subsection{The Difficulty}

We will again verify the criteria in the
Matching Lemma from \S \ref{MATCH}.  We
will change notation to reflect that we
have been calling $\widehat p/\widehat q$
the core predecessor of $p/q$

The general idea of the proof here is similar to
what we did in the previous chapter, but here
we must work with a weaker estimate, namely
\begin{equation}
\label{estimate0}
|P-\widehat P|=|Q-\widehat Q|=\frac{4\kappa}{\omega\widehat \omega}.
\end{equation}
The poorer quality of the estimate in Equation \ref{estimate0}
forces us to work harder in certain spots of the proof.

\subsection{Weak Horizontal Alignment}. 

 In our
context here, recall that
$\Upsilon$ is vertical translation by
$\omega/2-\widehat \omega/2$. 

\begin{lemma}
\label{unordered}
$\Upsilon(0,\widehat \tau)=(0,\tau)$
and $\Upsilon(0,\widehat \omega-\widehat \tau)=
(0,\omega-\tau)$.
\end{lemma}

\startproof
Since $\Upsilon$ is vertical translation
to $(\omega-\widehat \omega)/2$, we have
$$\Upsilon(\widehat \tau)=
\widehat \tau + (\omega-\widehat \omega)/2=
\tau.$$
The last equality is Statement 5 of
Lemma \ref{omnibus}.  
A similar calculation shows that
$\Upsilon(\widehat \omega-\widehat \tau)=\omega-\tau.$
Here are have abused notation and just shown
the action on the second coordinate.
\endproof

\begin{lemma}
$(\Sigma,\Pi)$ and $(\widehat \Sigma,\widehat \Pi)$
are weakly horizonally aligned.
\end{lemma}

\startproof
We take $\widehat H$ to be the line $y=\widehat \tau$
and we take $H$ to be the line $y=\tau$.
Since these lines have capacity $2$ w.r.t. the
relevant parameters, they only intersect
the plaid tilings in the middle of the first
interval.  The second intersection points
are outside $\widehat \Sigma$ and $\Sigma$
respectively.
\endproof

\subsection{Geometric Alignment}

\begin{lemma}
$(\Sigma,\Pi)$ and
$(\widehat \Sigma,\widehat \Pi)$ are
geometrically aligned.
\end{lemma}

\startproof
Let $\widehat z$ and $z$ be
corresponding vertical intersection
points.
Let $\widehat U$ and $U$ be the two
vertical intervals respectively containing
$\widehat z$ and $z$.  We need to prove that
$\Upsilon(\widehat U)=U$.

Let $\widehat L$ and $L$ be two slanting
lines of the same type which contain
$\widehat z$ and $z$ respectively and
have the same $y$-intercept.
The two lines $\Upsilon(\widehat L)$ and $L$
have the same $y$ intercept.  The difference
in their slopes is
\begin{equation}
\frac{4\kappa}{\omega \omega'}.
\end{equation}
Hence
\begin{equation}
\|z-\Upsilon(\widehat z)\| 
\leq \frac{4\kappa\widehat \tau}{\omega \widehat \omega}<\frac{2}{\widehat \omega}.
\end{equation}
Here is a derivation of the final inequality.
Statement 4 of Lemma \ref{omnibus} implies
that $\widehat \kappa=0$.  Hence
$\widehat \tau<\widehat \omega/3$.  But then
$$4\kappa \widehat \tau<\frac{4}{3}\kappa \widehat\omega<2\omega.$$
The last inequality comes from Statement 7 of
Lemma \ref{omnibus}.

The estimate is not quite good enough.  The problem
is that $\widehat z$ is exactly $m/\widehat \omega$
from the endpoint of $\widehat U$, but it
could happen that $m=1$.  However, given the
congruence properties of $\tau$, and the slopes
of the $\cal P$ and $\cal Q$ lines, we see that
$m=1$ if and only if $\widehat U$ is contained
in one of the two lines of capacity $2$.  
If $V$ is not one of these lines, we see that
$\widehat z$ lies at least $2/\widehat \omega$
from the endpoints of $\widehat U$.  Hence
$\Upsilon(\widehat z)$ lies at
least $2/\widehat \omega$ from the endpoints
of $\Upsilon(\widehat U)$.  But then our
estimate tells us that $z \in \Upsilon(\widehat U)$ as
well.  Hence $U=\Upsilon(\widehat U)$.

Now we deal with the case when $V$ is contains
the right edges of $\Sigma$ and $\widehat \Sigma$.
The equation for $V$ is $x=\widehat \tau$.
Lemma \ref{right} tells us that $V$ has capacity
$4\kappa+2$ w.r.t. $p/q$.   This means that
the $y$-coordinate of $z$ has the form
$(2\kappa+1)/\omega$ mod $\Z$.  This will give
us a good estimate provided that $\omega$ is large enough.
To see that $\omega$ is large enough for a good estimate, note
the formula
$$\omega=\omega'+2\kappa \widehat \omega.$$
Hence $\omega \geq 6 \kappa +3$, and
$z$ is at least $(2\kappa+1)/\omega$ from 
the endpoint of $U$.   If $\Upsilon(\widehat z)$
and $z$ were contained in different vertical
unit intervals, we would have
$$
\frac{4\kappa \widehat \tau}{\omega\widehat \omega} \geq
\frac{1}{\widehat \omega}+\frac{2\kappa+1}{\omega}.
$$
This is impossible, because $\widehat \tau/\widehat \omega<1/2$.
(We don't even need the $1/\widehat \omega$ term on the right.)
So, even in this special case, we see that
$\Upsilon(\widehat z)$ and $z$ are in the same
unit vertical interval. 
\endproof

\subsection{Alignment of the Capacity Sequences}
\label{signcap}

Define
\begin{equation}
\lambda=\widehat \omega/\omega, \hskip 40 pt
\Psi(i)=|2\widehat \omega Pi - 2\widehat \omega  \widehat P i|.
\end{equation}
As long as Equation \ref{key0} holds, namely
$$
\Psi(i)<\ell + 2\lambda, \hskip 30 pt
\min(| \widehat c_i|, \omega-|\widehat c_i|) \geq \ell
$$
the signs of $ \widehat c_i$ and $C_i$ are the same.

Using Equation \ref{estimate0} and the
fact that $i \in [1,\widehat \tau]$ we see
that
\begin{equation}
\label{bound5}
\Psi(i)=\frac{8\kappa i}{\omega} \leq \frac{8\kappa \widehat \tau}{\omega}.
\end{equation}

Statement 4 of Lemma \ref{omnibus} tells us
that $\widehat \kappa=0$, and this forces
$\widehat \tau \leq \widehat \omega/3$.
Combining this with Statement 7 of
Lemma \ref{omnibus} we get $\Psi(i)<4$.
The only cases we need to worry about are
\begin{itemize}
\item $\widehat C_i=\pm 2$,
\item $\widehat C_i=\pm(\widehat \omega-1)$.
\item $\widehat C_i=\pm(\widehat \omega-3)$.
\end{itemize}

\noindent
{\bf Case 1:\/}
We will consider the case when $\widehat C_i=2$.
The case when $\widehat C_i=-2$ has the same
treatment.
If $\widehat C_i=2$ then $i=\widehat \tau$ and
the corresponding $\cal V$ line is the right edge of
$\widehat \Sigma$.  But then from Lemma \ref{right} we
know that the capacity of this line w.r.t is
$4\kappa +2$.  So, we either have $C_i=4\kappa+2$
as desired or $C_i=-(4\kappa+2)$.  We will suppose
that $C_i=-(4\kappa+2)$ and derive a contradiction.
In this case, we would have
$$
\Psi(i)=
2+\lambda(4\kappa+2)=
\frac{2\omega+4\kappa \widehat \omega +2\widehat \omega}{\omega}.$$
Combining this with Equation \ref{bound5} we would get
$$2\omega+4\kappa \widehat \omega + 2 \widehat \omega<
8 \kappa \widehat \tau.$$
Statement 4 of Lemma \ref{omnibus} says that 
$\widehat \kappa=0$, which forces $\widehat \tau<\widehat \omega/3$,
and this contradict the equation above.
\newline
\newline
{\bf Case 2:\/}
Lemma \ref{anchor} tells us that the vertical lines of
capacity $\pm(\widehat \omega-1)$ either occur
at $x=\pm \widehat \tau/2$ or $x=\pm(\widehat \omega-\widehat \tau)/2$,
mod $\widehat \omega$. 
The latter case is irrelevant: The lines lie outside $\widehat \Sigma$.
In the former case, the line of interest to is $x=\widehat \tau/2$.
In other words $i=\widehat \tau/2$.  Note that $2i=\widehat \tau$.
The case $2i=\widehat \tau$ is the one just considered.

We will suppose that $\widehat C_{2i}=2$.  The case
$\widehat C_{2i}=-2$ has the same treatment.
When $\widehat C_{2i}=2$, it means 
that
$$[\widehat\omega(2\widehat P)(2i)]_{2\widehat \omega}=2.$$
But this is the same as saying that
$$\widehat \omega(2\widehat P)(2i)=2\theta\widehat \omega+2,$$
for some integer $\theta$.
But then
$$\widehat \omega(2\widehat Pi)=1+\theta\widehat \omega.$$
Since the capacities are all even, this is only
possible if $\theta$ is odd.  But then
$$\widehat C_i=[\widehat \omega(2\widehat Pi)]_{2\widehat \omega}=-\widehat \omega+1.$$
So, $\widehat C_i$ and $\widehat C_{2i}$ have
opposite signs.  A similar argument shows that
$C_i$ and $C_{2i}$ have opposite signs. Case 2 now
follows from Case 1.
\newline
\newline
{\bf Case 3:\/}
When $\widehat \tau$ is even, the relevant
line of capacity $\pm(\widehat \omega-3)$ is the
line $x=3\widehat \tau/2$.  
Since $3\widehat \tau/2<\widehat \omega/2$,
this line is the one which intesects the
first block and is closer to the $y$-axis.
But $3\widehat \tau/2>\widehat \tau$, and
so our line does not intersect $\widehat \Sigma$.

When $\widehat \tau$ is odd, the relevant
line of capacity $\pm(\widehat \omega-3)$ is
the line 
$$x=\widehat \omega/2 - 3\widehat \tau/2.$$
It would happen that this line intersects
$\widehat \Sigma$.  We will deal with this
line in a backhanded way.

We already know that $8\kappa \widehat \tau/\omega<4$.
If we really have a sign change in this case,
we must have $\Psi(i)>3$.
But this means that
\begin{equation}
\label{bound6}
\frac{8\kappa \widehat \tau}{\omega} \in (3,4).
\end{equation}
we have
\begin{equation}
\Psi(i)=
\frac{8\kappa i}{\omega} =
\frac{4\kappa \widehat \omega}{\omega}-\frac{12 \kappa \widehat \tau}{\omega}.
\end{equation}
By Statement 7 of Lemma \ref{omnibus}, the first term
on the right hand side of this equation is at most $6$.
By Equation \ref{bound6}, the second term on the right
lies in $(9/2,6)$.  Therefore
$\Psi(i) \leq 3/2$, and we have a contradiction.
Hence, there are no sign changes, and
the two capacity sequences are arithmetically aligned.

\subsection{Calculating some of the Masses}

As a prelude to checking that the
mass sequences are arithmetically aligned,
we take care of some special cases. This
will make the general argument easier.

We define the mass and sign of an integer point
on the $y$-axis to be the mass and sign of the
$\cal P$ and $\cal Q$ lines containing them.
Unless explicitly stated otherwise, these
quantities are taken with respect to the
parameter $p/q$.  In this section we
establish some technical results about some
of the masses and signs.

\begin{lemma}
\label{omni2}
Let $\widehat p/\widehat q$ be the
core predecessor of $p/q$.  Let
$\Upsilon$ be the translation from
the Core Copy Lemma.
\begin{enumerate}
\item Let $\Omega$ be the vertical interval
of length $\widehat \omega$ and
centered at the point $(0,\omega/2)$.
The only points of mass less than
$4\kappa+2$ contained in $\Omega$ are
the points of mass $1$.
\item $\Upsilon$ maps the points of
mass $1$ w.r.t $\widehat p/\widehat q$ 
to the points of mass $1$ w.r.t.
$p/q$, and in a sign-preserving way.
\item $\Upsilon$ maps the points of
mass $\widehat \omega-2$ w.r.t $\widehat p/\widehat q$ 
to the points of mass $\omega-2$ w.r.t.
$p/q$, and in a sign-preserving way.
\end{enumerate}
\end{lemma}

We will prove Lemma \ref{omni2} through a series
of smaller steps. 
We first dispense with a minor technical point.

\begin{lemma}
\label{triv}
$3\widehat \tau<\widehat \omega$.
\end{lemma}

\startproof
Since $\widehat p/\widehat q$ is assumed
to be the nontrivial core predecessor of
$p/q$, Statement 4 of Lemma \ref{omnibus}
gives $\widehat \kappa=0$.  This forces
$3\widehat \tau \leq \widehat \omega$.
The case of equality would force
$\widehat p/\widehat q=1/2$.  But
then the even predecessor of
$p/q$ would be the even predecessor of
$1/2$, by Statement 3 of Lemma \ref{omnibus}.
This is not possible.
\endproof

\noindent
{\bf Proof of Statement 1:\/}
The points $(0,\tau)$ and
$(0,\omega-\tau)$ are the
two points of mass $1$.
We will give the proof when the sign
of $(0,\tau)$ is positive.  The other
case is essentially the same.  

Since
$\Omega$ is symmetrically placed with
respect to the horizontal midline of the
first block $[0,\omega]^2$, it suffices
to show that $\Omega$ does not contain
any points of positive sign and mass
$3,5,7,...,4\kappa+1$.  

The endpoints of $\Omega$ are
$$(0,\omega/2-\widehat \omega/2), \hskip 30 pt
(0,\omega/2+\widehat \omega/2).$$
Our proof refers to the work in
\S \ref{hard}.  In particular,
$h=\omega-2\tau$ and $w=\widehat \tau$.
Since the masses of the points on the
Since the signed masses of the points
on the $y$-axis occur in an
arithmetic progression mod $\omega$,
the points on the $y$-axis having
positive sign and mass $2\lambda+1$,
at least for $\lambda=1,...,2\kappa$, are
\begin{equation}
(0,\tau-\lambda h)+(0,\omega \Z).
\end{equation}
The second summand is included just so that we remember
that the whole assignment of masses is invariant under
translation by $(0,\omega)$.  Also, we say ``at least''
because once $\lambda$ is large enough the sign will
change. So, it suffices to prove that
$$\tau-\lambda h \in
[-\omega/2+\widehat \omega/2,\omega/2-\widehat \omega/2],
\hskip 30 pt \lambda=1,...,2\kappa.$$
Since these points occur in linear order, it suffices to prove
the following two inequalities:
\begin{equation}
\label{MAIN2}
\tau-h<\omega/2-\widehat \omega/2, \hskip 30 pt
\tau-2\kappa h>-\omega/2+\widehat \omega/2.
\end{equation}

Since $h=\omega-2\tau$, the first identity
is equivalent to
$3\omega-6\tau>\widehat \omega.$
But, by Statement 5 of Lemma \ref{omnibus}, 
$$3\omega-6\tau=3\widehat \omega-6\widehat \tau.$$
But the inequality
$$3\widehat \omega-6\widehat \tau>\widehat \omega.$$
is the same as the one proved in Lemma \ref{triv}.
This takes care of the first inequality.

It follows from Equation \ref{MAIN} in
\S \ref{hard} that
$$\tau-2\kappa h=-\omega+\tau+h+2w.$$ Hence,
the second inequality is the same as
\begin{equation}
\label{final}
(-\omega+\tau+h+2w)-(-\omega/2+\widehat \omega/2)>0.
\end{equation}
Plugging in $h=\omega-2\tau$ and $w=\widehat \tau$,
and using the relation $\omega-2\tau=\widehat \omega-2\widehat \tau$,
we see that the left hand side of Equation
\ref{final} is just $\widehat \tau$.
\endproof

\noindent
{\bf Remark:\/} As a byproduct of our proof we note that
the points $(0,\tau)$ and $(0,\tau+h+2w)$ have the same sign.
Also, $h+2w=\widehat \omega$. 
\newline

\noindent
{\bf Proof of Statement 2:\/}
In view of Lemma \ref{unordered},
it suffices to prove that
$(0,\widehat \tau)$ has the
same sign w.r.t. $\widehat p/\widehat q$
as $(0,\tau)$ has w.r.t. $p/q$.
We will consider the case
when $(0,\widehat \tau)$ has
positive sign w.r.t $\widehat p/\widehat q$.
The opposite case has essentially
the same treatment.

The argument from Case 1 from
\S \ref{signcap} tells us that
the horizontal lines through
$\widehat y_+$ have capacity
$4\kappa+2$ w.r.t. $p/q$.
Since this is twice an odd integer,
Lemma \ref{anchor} tells us that the
slanting lines through
$\widehat y_+$ have mass $2\kappa+1$.
Now we observe the following.
\begin{equation}
\label{masschange}
P\omega'=P'\omega' + (P-P') \omega' = 2p'  \pm
\frac{2}{\omega} \equiv
\frac{2}{\omega}  \hskip 10 pt {\rm mod\/}\hskip 5 pt 2\Z.
\end{equation}
The calculation in Equation \ref{masschange}
implies that the masses w.r.t. $p/q$ along the
$y$-axis (so to speak) change by
$\pm 2$ when the $y$-coordinate
changes by $\omega'$.  Since
$$\tau=\widehat \tau + k\omega'$$
we see that
the slanting lines through $\tau$
either have mass $1$ or have mass
$$4\kappa+1 \hskip 10 pt
{\rm mod\/} \hskip 5 pt \omega.$$
But 
$$\omega=\widehat \omega+2\kappa \omega'>6\kappa+2.$$
Hence either the slanting lines through
$\tau$ have mass $1$ w.r.t $p/q$ or they have
mass $4\kappa+1$.  But we already
know from Lemma \ref{unordered}
that these lines either have mass
$1$ or $-1$.  Since
$4\kappa+1 \not = -1$ we know
that the slanting lines through
$\tau$ have mass $1$, as desired.
\endproof

\noindent
{\bf Proof of Statement 3:\/}
Let $z_{\pm}$ denote the points on
the $y$-axis such that the slanting
lines through these points have
mass $\pm (\omega-2)$.  Likewise
define $\widehat z_{\pm}$.
We will treat the case when $(0,\tau)$ has positive
sign. The other case has the same proof.

The function $f(x) \to [\widehat Px+1]_2$ is locally affine, when
the domain is interpreted to be the circle $\R/2\Z$.
In the case at hand, this function changes by
$+2$ when we move from $\widehat \tau$ to 
$\widehat \omega-\widehat \tau$.  Given this
property, and the fact that $f(0)=\pm \widehat \omega$,
we see that 
$$f(\widehat \omega-\widehat 2\tau)=-\widehat \omega+2.$$

But we also know that $y_+=\tau$ and
$y_-=\omega-\tau$ by the previous result.
So, the same argument gives
$$f(\omega- 2\tau)= - \omega+2.$$
This shows that
$m_i=-\omega+2$ when
$\widehat m_i=-\widehat \omega+2$.
By the same argument, or symmetry, we see that
$m_i=\omega-2$ when
$\widehat m_i=\widehat \omega-2$.
\endproof

\subsection{Alignment the Mass Sequences}

We proceed as in \S \ref{masschange1}.
Let $\widehat y=\widehat \omega/2$ and
$y=\omega/2$.
Even though $\widehat y$ is an integer point,
note that
\begin{equation}
[\widehat P \widehat y +1]_2=[\widehat p+1]_2,
\hskip 30 pt
[P y +1]_2=[p+1]_2.
\end{equation}
Since $p$ and $\widehat p$ are integers having
the same parity, the two expressions above are
the same.  Hence, there is an integer
$\theta$ so that
$$
\widehat \omega
\widehat P\widehat y = \widehat \omega Py + \theta \widehat \omega.
$$

Recall that $i_0=(\omega-\widehat \omega)/2$
and $y=\widehat y+i_0$.  
We can re-write the last equation as
$\Psi(\widehat \omega/2)=0$, where
\begin{equation}
\Xi(i)=|(\widehat \omega\widehat Pi -\widehat \omega P(i+i_0)-
\theta\widehat \omega|
\end{equation}
The function $\Xi$ has the same properties as the
similar function from Equation \ref{key1}:
As long as
$$
\Xi(i)<\ell + \lambda, \hskip 30 pt
\min(|\widehat M_i|, \widehat \omega-|\widehat M_i|) \geq \ell
$$
the signs of $\widehat M_i$ and $M_i$ are the same.

We say that the index $i$
is {\it central\/} if $i \in (0,\widehat \omega)$.
Otherwise, we call $i$ {\it peripheral\/}.

\begin{lemma}
If $i$ is a central index, and $m_i'$ is
term with a sign, then 
$\widehat m_i$ and $m_i$ have the same sign.
\end{lemma}

\startproof
When $i$ is a central index, we have 
$|i-\widehat \omega/2| \leq \widehat \omega/2$.
Combining this with Equation \ref{estimate0}, we get
\begin{equation}
\Psi(i) \leq \frac{4\kappa}{\widehat \omega \omega} \times
\frac{\widehat \omega}{2} = \frac{2\kappa \widehat \omega}{\omega}.
\end{equation}
Statement 7 of Lemma \ref{omnibus} says that
$\kappa \widehat \omega<(3/2)\omega$.  Combining
this with the previous equation,
we see that $\Psi(i)<3$.  This means that there are
no sign changes unless $\widehat m_i = \pm 1$
or $\widehat m_i = \pm(\omega-2)$.  But Statements 2
and 3 of Lemma \ref{omni2} take care of these
special cases.
\endproof

\begin{lemma}
If $i$ is a peripheral index then
$\widehat m_i$ and $m_i$ have the same sign.
\end{lemma}

\startproof
Note that all the peripheral indices are
signed terms, because none of these
indices are in $\widehat \omega \Z$.

The index $j$ of each peripheral term $m_j$ differs
from the index $i$ of some central term by
$\pm \widehat \omega$.  We set 
\begin{equation}
j=i+\widehat \omega.
\end{equation}
The case when $j=i-\widehat \omega$ has
a similar treatment.

Consider the situation w.r.t. the
parameter $\widehat p/\widehat q$ first.
These terms repeat every $\widehat \omega$.
This meams that $m'_j=m'_i$.

Now observe that
\begin{equation}
\label{gapX}
[P\widehat \omega]_2=
[\widehat P \widehat \omega + (P-\widehat P)\omega]_2=
[(P-\widehat P)\omega]_2 =
\pm \frac{4\kappa}{\omega}.
\end{equation}
Therefore
\begin{equation}
m_j=m_i \pm 4\kappa.
\end{equation}

By Statement 1 of Lemma \ref{omni2}, the central
terms w.r.t. $p/q$ 
either equal $\pm 1$ or are greater than
$4\kappa+2$.  So, if $m_i \not = \pm 1$,
we have the terms $m_i$ and $m_j$ have
the same sign.  But then
$$\sigma(m_j)=\sigma(m_i)=\sigma(m'_i)=\sigma(m'_j).$$
Here $\sigma$ denotes the sign. The
middle equality comes from the previous result
about the central indices.

What about when $m_i'=\pm 1$?
In this case (since $j=i+\widehat \omega$)
we must have $i=\tau$ and $j=\tau+\widehat \omega$.
The remark following the proof of Statement 1
of Lemma \ref{omni2}.
takes care of this case.
\endproof

Now we know that $(\Sigma,\Pi)$
and $(\widehat \Sigma,\widehat \Pi)$ are
arithmetically aligned.  We have
verified the conditions of the
Matching Lemma, and so
$\Sigma \cap \Pi=\Upsilon(\widehat \Sigma \cap \widehat \Pi)$.
This completes the proof of the
Core Copy Lemma.

\newpage

\section{References}

[{\bf DeB\/}] N. E. J. De Bruijn, {\it Algebraic theory of Penrose's nonperiodic tilings\/},
Nederl. Akad. Wentensch. Proc. {\bf 84\/}:39--66 (1981).
\newline
\newline
[{\bf H\/}] W. Hooper, {\it Renormalization of Polygon Exchange Transformations
arising from Corner Percolation\/}, Invent. Math. {\bf 191.2\/} (2013) pp 255-320
\newline
\newline
[{\bf S0\/}] R. E. Schwartz, {\it Inroducing the Plaid Model\/},
preprint 2015.
\newline
\newline
[{\bf S1\/}] R. E. Schwartz, {\it Outer Billiard on Kites\/},
Annals of Math Studies {\bf 171\/} (2009)

\end{document}